\documentclass[final,leqno]{siamltex}
\usepackage{threeparttable}
\usepackage{dcolumn}
\newcolumntype{d}{D{.}{.}{-1}}
\usepackage{subfigure}
\usepackage[german,english]{babel}
\usepackage{amssymb}
\usepackage{algorithm}
\usepackage{algorithmic}
\usepackage{color}
\usepackage{graphicx}
\usepackage{booktabs}
\usepackage{amsfonts}
\usepackage{soul}
\usepackage{mathrsfs}
\usepackage{mathtools}
\usepackage{amsmath, amssymb}
\usepackage{float}
\usepackage{threeparttable}
\usepackage{dcolumn}
\newcolumntype{d}{D{.}{.}{-1}}
\newtheorem{remark}[theorem]{Remark}
\usepackage{enumitem}
\usepackage[small,nohug,heads=vee]{diagrams}
\diagramstyle[labelstyle=\scriptstyle]

\newcommand*{\Scale}[2][4]{\scalebox{#1}{$#2$}}%

\title{Symplectic Model Reduction of Hamiltonian Systems\thanks{
This research was supported by the Office of Naval Research.}}

\author{Liqian Peng\thanks{Department of Mechanical and
Aerospace Engineering, and Institute for Networked Autonomous
Systems, University of Florida, Gainesville, FL 32611-6250
(liqianpeng@ufl.edu).}
        \and Kamran Mohseni\thanks{Department of Mechanical and
Aerospace Engineering, Department of Electrical and Computer
Engineering, and Institute for Networked Autonomous Systems,
University of Florida, Gainesville, FL 32611-6250
(mohseni@ufl.edu).}}

\begin{document}

\maketitle
\begin{abstract}
In this paper, a symplectic model reduction technique, proper symplectic decomposition (PSD) with symplectic Galerkin projection, is proposed to save the computational cost for the simplification of large-scale Hamiltonian systems while preserving the symplectic structure.  As an analogy to the classical proper orthogonal decomposition (POD)-Galerkin approach, PSD is designed to build a symplectic subspace to fit empirical data, while the symplectic Galerkin projection constructs a reduced Hamiltonian system on the symplectic subspace. For practical use, we introduce three algorithms for PSD, which are based upon:  the cotangent lift,   complex singular value decomposition, and   nonlinear programming.  The proposed technique has been proven to preserve system energy and stability.  Moreover, PSD can be combined with the discrete empirical interpolation method to reduce the computational cost for nonlinear Hamiltonian systems.  Owing to these properties, the proposed technique is better suited than the classical POD-Galerkin approach for model reduction of Hamiltonian systems, especially when long-time integration is required.  The stability, accuracy, and efficiency of the proposed technique are illustrated through numerical simulations of linear and nonlinear wave equations.
\end{abstract}



\pagestyle{myheadings} \thispagestyle{plain} \markboth{LIQIAN PENG
AND KAMRAN MOHSENI}{SYMPLECTIC  MODEL REDUCTION OF HAMILTONIAN SYSTEMS}

\section{Introduction}
To save computational costs, model reduction seeks to approximate high-dimensional dynamical systems using simpler, lower order ones that can capture the dominant dynamic properties.  The need for model reduction arises because, in many cases, direct numerical simulations are often so computationally intensive that they either cannot be performed as often as needed or are only performed in special circumstances.  See~\cite{SorensenDC:01a} for a survey on the available model reduction techniques.

Among these techniques, the proper orthogonal decomposition (POD) with Galerkin projection, which was first introduced by Moore~\cite{MooreBC:81a}, has wide applications in many fields of  science and engineering, such as electric circuit analysis~\cite{Marsden:99e}, structural dynamics~\cite{AmabiliM:03a}, and fluid mechanics~\cite{HolmesP:02a,RowleyCW:04a}, to list a few.  As an empirical model reduction technique, the POD-Galerkin approach (or POD for short) involves an offline-online splitting methodology.  In the offline stage, empirical data is generated by direct numerical simulations of the original system. If the original system is represented by a PDE, a discretized high-dimensional model can be derived by the finite difference, finite element, and finite volume methods. The POD can be applied to compute an optimal subspace to fit the empirical data.  A reduced system is then constructed by projecting the high-dimensional system to this subspace.  In the online stage, one can solve the reduced system in the low-dimensional subspace.  Recently, many variants of POD-Galerkin have been developed to reduce the complexity of evaluating the nonlinear term of the vector field, such as  trajectory piecewise linear approximation~\cite{WhiteJ:03a},  missing point estimation~\cite{WillcoxK:08a},  Gappy POD~\cite{WillcoxK:04a,WillcoxK:06a}, empirical interpolation method~\cite{PateraAT:04a,PateraAT:06a},  and discrete empirical interpolation method (DEIM)~\cite{SorensenDC:10a,DrohmannM:12a}. Thanks to these methods, the computational complexity during the online stage is independent of the dimension of the high-dimensional model.

Generally, the classical POD is not guaranteed to yield a stable reduced system, even if the original system is stable~\cite{PetzoldLR:03a,PrajnaS:03a}. The instability of a reduced system is often accompanied by blowup of system energy and flow volume. Therefore, when the original large-scale system is conservative, it is preferable to construct a low-dimensional reduced system that preserves the geometric structure and allows symplectic integrators.  However, much less effort has been expended in the field of geometric model reduction. In the context of a Lagrangian system, Lall et al.~\cite{Marsden:03h} showed that performing a Galerkin projection on the Euler--Lagrange equation and lifting it to the tangent bundle of the phase space leads to a reduced system that preserves Lagrangian structure. In order to reduce the complexity of nonlinear Lagrangian systems, Carlberg et al. combined Lall's method with the Gappy POD to derive reduced nonlinear Lagrangian systems~\cite{CarlbergK:14a}. In the control community, the balanced truncation~\cite{HartmannC:10a}, moment matching~\cite{vdSchaft:10a}, and tangential interpolation~\cite{vdSchaft:12a} approaches were used to preserve the port-Hamiltonian structure.

In this paper, we propose a new model reduction technique, proper symplectic decomposition (PSD), that preserves the symplectic structure underlying the Hamiltonian mechanics. Our main focus is to develop a basic framework behind symplectic model reduction, which allows us to derive energy preservation and stability preservation. The proposed technique yields reduced Hamiltonian systems which are applicable to long-time integration.   Compared with other empirical model reduction algorithms that preserve system energy, the PSD is easier for applications; the computation complexity can be the same magnitude as the original POD and DEIM for both offline and online stages.  The PSD also increases the flexibility to construct an optimal subspace that can yield a more accurate solution for the same subspace dimension.

The remainder of this paper is organized as follows. Preliminaries of Hamiltonian systems and symplectic integrators are briefly reviewed in Section \ref{sec:background}. Section \ref{sec:projection} presents the symplectic projection, which constructs reduced Hamiltonian systems. In Section \ref{sec:pod}, three different algorithms are proposed to construct a symplectic matrix, including the cotangent lift, complex SVD, and nonlinear programming. In Section \ref{sec:deim}, the symplectic discrete empirical interpolation method (SDEIM) is developed in order to reduce the complexity of evaluating the nonlinear vector term. Sections \ref{sec:projection}, \ref{sec:pod}, and \ref{sec:deim} respectively associate with the classical Galerkin projection, POD, and DEIM. In Section \ref{sec:numerical}, the stability, accuracy, and efficiency of the proposed technique are illustrated through numerical simulations of linear and nonlinear wave equations.  Finally, conclusions are offered in Section \ref{sec:conclusion}.

\section{Hamiltonian System and Symplectic Integrator}\label{sec:background}
Let $\mathbb{V}$ be a vector space of dimension $2n$. A symplectic form on $\mathbb{V}$ is a closed, nondegenerate, skew-symmetric bilinear form,  $\Omega: \mathbb{V}\times \mathbb{V} \to \mathbb{R}$. The pair $(\mathbb{V}, \Omega)$ is called a symplectic vector space. Assigning a symplectic form $\Omega$ to $\mathbb{V}$ is referred to as giving $\mathbb{V}$ a symplectic structure. With canonical coordinates on $\mathbb{V}$ denoted by $(q_1, \ldots, q_n, p_1, \ldots,p_n)$,  $\Omega$ takes a canonical form, $\Omega = \sum\nolimits_{i = 1}^n {d{q_i} \wedge d{p_i}}$. Throughout this paper, we implicitly assume that $\mathbb{V}$ is defined over the field $\mathbb{R}$, which means $\mathbb{V}=\mathbb{R}^{2n}$. Moreover, for all  $v_1,v_2 \in \mathbb{V}$, $\Omega$ is represented by the Poisson matrix $J_{2n}$, i.e.,
\begin{equation*}
\Omega(v_1, v_2)=v_1^T J_{2n} v_2, \qquad {J_{2n}} = \begin{bmatrix}
  0_n & {I_n} \\
 {-I_n} & 0_n \\
 \end{bmatrix},
\end{equation*}
where $I_n$  is the $n\times n$  identity matrix. It is easy to verify that $J_{2n}J_{2n}^T=J_{2n}^TJ_{2n}=I_{2n}$, and $J_{2n}J_{2n}=J_{2n}^TJ_{2n}^T=-I_{2n}$, where the superscript $T$ represents the transpose of a matrix.

Let $H: \mathbb{V} \to \mathbb{R}$ denote a smooth Hamiltonian function. The time evolution of an autonomous Hamiltonian system is defined by:
\begin{equation}\label{hamiltons}
\dot q=  \nabla_p  H(q,p), \qquad
\dot p= -\nabla_q  H(q,p),
 \end{equation}
where $q=[q_1; \ldots; q_n]\in \mathbb{R}^{n}$, and $p=[p_1; \ldots; p_n]\in \mathbb{R}^{n}$. We abstract this formulation by introducing the phase space variable $x=[q;p]$\footnote{The notations $[q, p]$ and $[q; p]$ are the same as the corresponding functions in MATLAB.} and the abstract Hamiltonian differential equation
\begin{equation}\label{fom}
\dot x= X_H(x),
\end{equation}
where $X_H(x):=J_{2n} \nabla_x H(x)$ is the Hamiltonian vector field.  The flow $\Psi_t$ of $X_H$ is a symplectomorphism, meaning that it conserves the symplectic two-form $\Omega$, the system Hamiltonian $H$, and the volume of flow $\Theta$.  We recommend that readers refer to~\cite{Marsden:03j} for more details of these fundamental facts of symplectic geometry.

 Symplectic integrators are numerical schemes for solving a Hamiltonian system, while preserving the underlying symplectic structure. If the symplectic structure is preserved, then the flow volume and system energy are automatically conserved during time integration.   By virtue of these advantages, symplectic integrators have been widely applied to  long-time integrations of  molecular dynamics, discrete element methods,  accelerator physics, and  celestial mechanics~\cite{HairerE:06a}.

 Let $\delta t$ denote the unit step for time integration. The symplectic Euler methods
\begin{align*}
 q^{j+1}=q^{j}+\delta t \nabla_p
 H(q^{j+1}, p^j) \qquad  or \qquad
 q^{j+1}=q^{j}+\delta t \nabla_p
 H(q^{j}, p^{j+1})
 \end{align*}
 \vspace{-12mm}
\begin{align*}
 p^{j+1}=p^{j}-\delta t \nabla_q
 H(q^{j+1}, p^j)  \qquad \ \ \ \qquad
 p^{j+1}=p^{j}-\delta t \nabla_q
 H(q^{j}, p^{j+1})
 \end{align*}
 are symplectic integrators of order one.  They are implicit for general Hamiltonian systems. For separable $H(q,p)=T(p)+U(q)$, however, both variants turn out to be explicit. If the implicit midpoint rule is applied, then a second order symplectic scheme is obtained:
 \begin{equation}\label{sischeme}
{x^{j + 1}} = {x^j} + \delta t{J_{2n}}{\nabla _x}H\left( {\frac{{{x^{j + 1}} + {x^j}}}{2}} \right).
\end{equation}

\vspace{-2mm}
 Most of the usual numerical methods, such as the primitive Euler scheme and the classical Runge--Kutta scheme, are not symplectic integrators.  A comprehensive review of symplectic integrators and their applications for Hamiltonian ODEs can be found in~\cite{HairerE:06a, McLachlanRI:06a}; the extension for Hamiltonian PDEs can be found in~\cite{BridgesTJ:06a}, where some structure-preserving discretization methods are discussed to transform Hamiltonian PDEs into Hamiltonian ODEs.

\section{Symplectic Projection}\label{sec:projection}
The symplectic projection takes advantage of empirical data to construct a reduced system, while simultaneously preserving the underlying symplectic structure. In other words, if the original system is Hamiltonian, the reduced system remains Hamiltonian, but with significantly fewer dimensions.

\subsection{Definitions of Symplectic Projection}\label{sec:def} Let $(\mathbb{V}, \Omega)$ and $(\mathbb{W}, \omega)$ be two symplectic vector spaces; $\dim(\mathbb{V})=2n$, $\dim(\mathbb{W})=2k$, and $k \le n$.

\medskip
\begin{definition}
A {\textbf{symplectic lift}} is a linear mapping $\sigma: \mathbb{W} \to \mathbb{V}$ that preserves the symplectic structure:
 \begin{equation} \label{symp}
 \omega( z_1, z_2)=\Omega(\sigma (z_1), \sigma (z_2)),
 \end{equation}
for all  $z_1,z_2 \in \mathbb{W}$.
\end{definition}
\medskip

 Let $z \in \mathbb{W}$ and $x \in \mathbb{V}$. Using canonical coordinates, a symplectic lift $\sigma: z\mapsto x$ can be written as
 \begin{equation*}
 x=Az,
 \end{equation*}
 where  $A\in \mathbb{R}^{2n\times 2k}$ and satisfies
\begin{equation}\label{AJA}
A^T J_{2n} A =J_{2k}.
\end{equation} 
A matrix  that satisfies (\ref{AJA}) for some $k$ and $n$ with $k\le n$ is called a \emph{symplectic} matrix.  The set of all ${2n\times 2k}$ symplectic matrices is  the \emph{symplectic Stiefel manifold}, denoted by $Sp({2k},\mathbb{R}^{2n})$.
Moreover, since $J_{2n}$ and $J_{2k}$ are nonsingular, (\ref{AJA}) itself requires that $k\le n$ and ${\rm{rank}}({A})= 2k$.

\medskip
\begin{definition}
The {\textbf{symplectic inverse}} of a real matrix $A\in \mathbb{R}^{2n\times 2k}$, denoted as $A^+$, is defined by \begin{equation}\label{A+}
A^+=J_{2k}^TA^TJ_{2n}.
\end{equation}
\end{definition}
Although $A^+$ is not equal to the Moore-Penrose pseudoinverse $(A^TA)^{-1}A^T$ in general, $A^+$ has several interesting properties, as stated in the following two lemmas. Using the definition of $A^+$, it is straightforward to verify Lemma \ref{lem:sym1}.
  \medskip
 \begin{lemma}\label{lem:sym1}
Suppose $A\in \mathbb{R}^{2n\times 2k}$ and $A^+$ is the symplectic inverse of $A$. Then,
\begin{equation}\label{A++}
A=(A^+)^+,
\end{equation}
\vspace{-10mm}
\begin{equation}\label{A+T+T}
A=(((A^+)^T)^+)^T,
\end{equation}
\vspace{-10mm}
 \begin{equation}\label{exch}
A^+J_{2n}=J_{2k}A^T.
\end{equation}
\end{lemma}

\vspace{-4mm}

\begin{lemma}\label{lem:sym2}
Suppose $A\in \mathbb{R}^{2n\times 2k}$ and $A^+$ is the symplectic inverse of $A$. Then the following are equivalent:
\begin{enumerate}[label={\rm{(}}\alph*{\rm{)}}]
\setlength{\itemsep}{2pt}
\item $A\in Sp({2k},\mathbb{R}^{2n})$.
  \item $(A^+)^T\in Sp({2k},\mathbb{R}^{2n})$.
  \item $A^+A=I_{2k}$.
\end{enumerate}
\end{lemma}
 \medskip
 \begin{proof}
 $(a)  \Rightarrow (b):$    Replacing $A^+$ by (\ref{A+}) and using $A^TJ_{2n}A=J_{2k}$ yield
  \begin{equation*}
{A^ + }{J_{2n}}{({A^ + })^T}=(J_{2k}^TA^TJ_{2n})J_{2n}(J_{2n}^TAJ_{2k})=J_{2k}^T(A^TJ_{2n}A)J_{2k}=J_{2k}.
 \end{equation*}
Since $((A^+)^T)^T=A^+$, we have $(A^+)^T\in Sp({2k},\mathbb{R}^{2n})$.

$(b)  \Rightarrow (c):$ Since $(A^+)^T\in Sp({2k},\mathbb{R}^{2n})$, we have ${A^ + }{J_{2n}}{({A^ + })^T}=J_{2k}$. Substituting  $A$ by (\ref{A+T+T}) and simplifying the expression yeild
\begin{align*}
 {A^ + }A  &= {A^ + }{({({({A^ + })^T})^ + })^T} = {A^ + }{(J_{2k}^T{({({A^ + })^T})^T J_{2n} })^T} = {A^ + }{(J_{2k}^T{A^ + }{J_{2n}})^T} \\
  &  = {A^ + }J_{2n}^T{({A^ + })^T}{J_{2k}} =  - ({A^ + }{J_{2n}}{({A^ + })^T}){J_{2k}}=- {J_{2k}}{J_{2k}} = {I_{2k}}.
 \end{align*}

$(c)  \Rightarrow (a):$ Replacing $A^+$ by (\ref{A+}) and plugging it into $A^+A= {I_{2k}}$, we obtain $J_{2k}^TA^TJ_{2n}A=I_{2k}$. Left multiplying $J_{2k}$ on both sides of this equation yields $A^TJ_{2n}A=J_{2k}.$
 \hfill
\end{proof}
  \medskip
\begin{definition}
Let $z\in \mathbb{W}$ and $x\in \mathbb{V}$.
Using the canonical coordinates, a linear mapping $\pi: x  \mapsto z$ is a {\textbf{symplectic projection}} if there exists a symplectic matrix $A\in Sp({2k},\mathbb{R}^{2n})$, such that
 \begin{equation}\label{projection}
z=A^+x.
\end{equation}
\end{definition}
\vspace{-5mm}

\begin{remark}
If we generalize $(\mathbb{W}, \Omega)$ and $(\mathbb{V}, \omega)$ to two symplectic manifolds and consider nonlinear transformations,  the symplectic lift and symplectic projection respectively correspond to the symplectic embedding and symplectic submersion in symplectic geometry. Since this article focuses on providing efficient numerical algorithms for practical applications, we only consider linear transformations between two vector spaces, although both the original and reduced systems can be nonlinear.
\end{remark}
\medskip

By Lemma \ref{lem:sym2}, $\pi \circ \sigma$ is the identity map on $\mathbb{W}$. Now suppose $x=Az$, where $A^TJ_{2n}A=J_{2k}$.
 Using the chain rule, we obtain ${\nabla _z}H(Az)=A^T{\nabla _x}H(x)$.
 Taking the time derivative of (\ref{projection}) and using (\ref{fom}) and (\ref{exch}), the time evolution of $z$ is given by
\begin{equation*}
\dot z = {A^ + }\dot x = {A^ + }{J_{2n}}{\nabla _x}H(x) = J_{2k} A^T{\nabla _x}H(x)= J_{2k}{\nabla _z}H(Az),
\end{equation*}
where the last expression is a Hamiltonian vector field.

 \medskip
\begin{definition}
The {\textbf{symplectic Galerkin projection}}, or {\textbf{symplectic projection}}, of a $2n$-dimensional Hamiltonian system $\dot x=J_{2n}\nabla_x H(x),$ with an initial condition  $x(0)=x_0$  is given by a $2k$-dimensional ($k\le n$) system
 \begin{equation}\label{rom}
\dot z = J_{2k}{\nabla _z}\tilde H(z); \qquad z_0=A^+x_0,
\end{equation}
where $\tilde H:=H \circ A$ is the reduced Hamiltonian function, $A\in Sp({2k},\mathbb{R}^{2n})$ is a symplectic matrix, and $A^+=J_{2k}^TA^TJ_{2n}$ is the symplectic inverse of $A$.
\end{definition}
 \medskip

 \begin{remark}
 Some Hamiltonian systems, such as the Burgers equation and KdV equation, have nontrivial symplectic structures~\cite{Marsden:03j}, which can be written in the form
\begin{equation}\label{gfom}
\dot x= J \nabla_x H(x),
\end{equation}
where $J\in \mathbb{R}^{2n\times 2n}$ is a nondegenerate skew-symmetric matrix. Especially, when $J =J_{2n}$, (\ref{gfom}) denotes a standard Hamiltonian system.  Notice that for any nondegenerate skew-symmetric matrix $J$, there exists a congruent transformation such that $J=Q^T J_{2n} Q$~\cite{EvesH:80a}, where $Q$ is a non-singular matrix. Let $y=(Q^T)^{-1}x$, then $\nabla_x H(x)=Q^{-1}\nabla_y H(Q^T y)$. It follows that
$$\dot y= (Q^T)^{-1}  \dot x= (Q^T)^{-1} {J} \nabla_x H(x)= J_{2n}\nabla_y H(Q^T y).$$
The last equation indicates that a Hamiltonian equation with nontrivial symplectic structures (\ref{gfom}) can be transformed to a standard Hamiltonian equation, and therefore can be simplified by the symplectic projection.
 \end{remark}

\subsection{Linear Hamiltonian Systems}
A Hamiltonian system is \emph{linear} if $H(x) = {\textstyle{1 \over 2}}{x^T}Lx$,  where  $L$ is a  ${2n\times 2n}$ real symmetric matrix. Let $K:=J_{2n}L$, the linear Hamiltonian system can be written as
\begin{equation}\label{LTI}
\dot x=J_{2n}Lx=Kx.
\end{equation}
A matrix of the form $K=J_{2n}L$, where $L$ is symmetric, is called a \emph{Hamiltonian matrix}. In addition, the set of all $2n\times 2n$ Hamiltonian matrices, denoted by $\mathfrak{sp}(\mathbb{R}^{2n})$, is a Lie algebra~\cite{MeissJD:07a}. The fundamental matrix solution to (\ref{LTI}) is given by
\begin{equation}
x(t)=e^{Kt}x_0.
\end{equation}
Since $\exp(Kt)$ satisfies $(\exp(Kt))^TJ_{2n}\exp(Kt)=J_{2n}$, we have $\exp(tK)\in Sp(2n, \mathbb{R}^{2n})$, which means that the matrix exponential of a Hamiltonian matrix is symplectic. Conversely, the logarithm of a square symplectic matrix is Hamiltonian; see reference~\cite{MeyerKR:09a} for the proof.

Applying the symplectic projection on the linear system (\ref{LTI}) gives
\begin{equation}\label{rLTI}
\dot z = {A^ + }\dot x = {A^ + }Kx = A^+KAz=\tilde Kz,
\end{equation}
where $\tilde K:=A^+KA$ is the reduced linear operator.  Using (\ref{exch}) in Lemma \ref{lem:sym1} gives
\begin{equation}\label{symsym}
\tilde K=A^+(J_{2n}L)A=J_{2k}(A^TLA)=J_{2k}\tilde L.
\end{equation}
Since $\tilde L=A^TLA$ is symmetric, we have $\tilde K\in \mathfrak{sp}(\mathbb{R}^{2k})$, which implies that the reduced linear system (\ref{rLTI}) is also Hamiltonian.
 
 Furthermore, a reduced system can also be obtained by directly plugging the reduced Hamiltonian $\tilde H(z)$ into (\ref{rom}); for the linear case, $\tilde H(z)={\textstyle{1 \over 2}}{(Az)^T}L(Az)$. Since the reduced system constructed by the symplectic projection is always Hamiltonian, consequently,  energy and stability are preserved during the time evolution.

\subsection{Energy Preservation}  Let $\Delta H(t):=H(x(t))-\tilde H(z(t))$ denote the energy discrepancy  between the state $x(t)$ and its approximation, $A z(t)$, derived from a reduced system. Since both the original and reduced systems are Hamiltonian, the system energy is conserved during time evolution. Moreover, $\tilde H=H \circ A$ by the definition. Thus,  $\Delta H(t)$ is  determined by the initial condition $x_0$ and the basis matrix $A$ for all $t$, i.e.,
  \begin{equation}\label{engerror}
  \Delta H(t)=H(x(0))-\tilde H(z(0))=H (x_0)-H (AA^+x_0).
\end{equation}
If $x_0\in {\rm {Range}}(A)$, we have $AA^+ x_0= x_0$, which implies $\Delta H(t)=0$ for all $t$; we say that the reduced system is \emph{energy preserving}.

If $x_0 \notin {\rm {Range}}(A)$,  we can always extend $A$ to a larger symplectic matrix $A_{\rm{ext}}$ such that  the reduced system remains energy preserving.  Specifically, suppose $A=[A_1, A_2]$ for $A_1, A_2 \in \mathbb{R}^{2n \times k}$. Since $x_0 \notin {\rm {Range}}(A)$,  we must have $r_0:=x_0-AA^+x_0\ne 0$.  Thus, the unit vector, $\hat r_0:=  r_0/\|r_0\|$, is well-defined. One possible extension of $A$ is given by
  \begin{equation} \label{barA}
  A_{\rm{ext}}=[A_1, \hat r_0, A_2, J_{2n}^T \hat r_0].
   \end{equation}
It is straightforward to verify that $A_{\rm{ext}}^T J_{2n}A_{\rm{ext}}=J_{2k+2}$, and $x_0-A_{\rm{ext}}A_{\rm{ext}}^+x_0=0$. The last equation means that $x_0 \in  {\rm {Range}}(A_{\rm{ext}})$, and therefore, $\Delta H(t)=H (x_0)- H (A_{\rm{ext}} A_{\rm{ext}}^+x_0)=0$ for all $t$.

\subsection{Stability Preservation}  Let $S$ denote an open set of $\mathbb{V}$ that contains $x_0$, and let $\partial S$ denote the boundary of $S$. Moreover, we assume $x_0\in {\rm {Range}}(A)$, and the initial condition of the reduced system is given by $z_0=A^+x_0$. The following two theorems imply that energy preservation is a strong indicator for preserving stability.

\medskip
\begin{theorem} \label{thm:bound}
Consider the Hamiltonian system (\ref{fom}) with the initial condition $x_0\in \mathbb{V}$. If there exists a bounded neighborhood $S$ of $x_0$ in $\mathbb{V}$ such that $H(x_0) < H(x)$, or $H(x_0)> H(x)$, for all $x\in \partial S$, then both the original system and the reduced system  constructed by the symplectic projection are bounded for all $t\in \mathbb{R}$.
\end{theorem}
\medskip

\begin{proof}
We first prove that the statement is true for the case that $H(x_0) < H(x)$ for all $x\in \partial S$.
Let $E=\min\{H(x):x\in \partial S\}$; so $H(x_0)<E$. Because of energy conservation, we have $H(x(t))=H(x_0)<E$ for all $t$. It follows that $x(t)\in S$ for all $t$, because if not, there is a time $t_1$ when $x(t_1)\in \partial S$, and $H(x(t_1))\ge E$, a contradiction.

Let $S_A=S \cap {\rm{Range}}(A)$. Since $S$ is a bounded open set in $\mathbb{V}$, $S_A$ is also open in ${\rm{Range}}(A)$ and bounded. Moreover, $\partial S_A= \partial S \cap {\rm{Range}}(A)$. Thus, $H(x_0)<H (x)$ for all $x\in S_A$. By the same argument in the last paragraph, we must have $A z(t)\in S_A$ for all $t$, which means that the reduced system constructed by the symplectic projection is also bounded.

Finally, if $H(x_0) < H(x)$ is replaced by $H(x_0)> H(x)$, we can define $\hat H(x)=-H(x)$. Then, we have $\hat H(x_0)< \hat H(x)$ for all $x \in \partial S$. By the same argument, the conclusion still holds.
\hfill
\end{proof}
\medskip

 When the original Hamiltonian system is linear, then the reduced system is also linear. The boundedness of linear systems implies that both the original and reduced systems are stable.

 \medskip
 \begin{theorem}\label{minimal}
If $x_*\in {\rm{Range}}(A)$ is a strict local minimum or maximum of $H$, then $x_*$ is a stable equilibrium for both the original Hamiltonian system and the  reduced system constructed by the symplectic projection.
 \end{theorem}
 \medskip

 \begin{proof}
   The Dirichlet's stability theorem states that if $x_*$ is a strict local minimum (or maximum) of $H$, then $x_*$ is a stable equilibrium for (\ref{fom})~\cite{MeyerKR:09a}. Suppose $S$ is a neighborhood of $x_*$, and $x_*$   is the minimum (or maximum) of $H$ in $S$. It immediately follows that $x_*$   is also the minimum (or maximum) of $H$ in $S_A$, where $S_A= S \cap {\rm{Range}}(A)$. Thus, by the Dirichlet's stability theorem, $x_*$ is also the stable equilibrium of the reduced Hamiltonian system.
   \hfill
 \end{proof}
 \medskip

The symplectic projection is analogous to the Galerkin projection, both of which construct reduced equations in some low dimensional subspaces. However, the symplectic projection yields a reduced symplectic system by (\ref{rom}) while the Galerkin projection generally destroys the symplectic structure. Evolving the system (\ref{rom}) by a symplectic integrator preserves system energy and stability. By contrast, even if the POD subspace can provide an accurate representation of the empirical data, the reduced system constructed by the Galerkin projection may not be able to preserve these properties of the dynamics.  In the next subsection, we shall discuss some PSD algorithms to construct a symplectic matrix $A$. This approach is an analogy to the POD that constructs an orthonormal basis matrix.

\section{Proper Symplectic Decomposition (PSD)}\label{sec:pod}
Let  $x(t_i)=[q(t_i);p(t_i)]\in \mathbb{R}^{2n}$ ($i=1,\ldots, N$) denote $N$  data points.   Define a snapshot matrix
\begin{equation}\label{snap}
M_x:=[x(t_1),\ldots,x(t_N)].
 \end{equation}
The symplectic projection of $M_x$ onto a low dimensional subspace is given by $M_z=A^+M_x$, where $A\in Sp({2k},\mathbb{R}^{2n})$, $M_z=[z(t_1),\ldots, z(t_N)]\in \mathbb{R}^{2k\times N}$, and  $z(t_i)= A^+ x(t_i)$. The same projection of $M_x$ in the original coordinates is given by
$AM_z$, or $AA^+M_x$.

The Frobenius norm $\| \cdot \|_F$ can be used to measure the error between $M_x$ and its projection $\tilde M_x$. Suppose a symplectic matrix $A$   minimizes the projection error in a least squares sense. Then, $A$ is a solution of the following  optimization problem:
\begin{equation}\label{opt}
\begin{aligned}
&{\rm{minimize}} \quad \ \|M_x-AA^+ M_x\|_F \\
 &{\rm subject \ to}  \quad  A^TJ_{2n}A=J_{2k}.
\end{aligned}
\end{equation}
 Since the objective function has a fourth-order term in $A$ after an expansion, (\ref{opt}) can only be solved iteratively.   Because matrix $A$ has $4nk$ elements, direct solving (\ref{opt}) is very expensive if $n \gg 1$.
For this reason, we propose three efficient algorithms to construct an approximated optimal solution for the symplectic matrix $A$: these are the cotangent lift,   complex SVD, and   nonlinear programming (NLP).

\subsection{Cotangent Lift}\label{sec:clm}
In this section, an SVD-based algorithm is proposed to construct a symplectic matrix directly. The idea is to search the optimal matrix, $A_1$, in a subset of $Sp(2k,\mathbb{R}^{2n})$, such that all the empirical data approximately lies on the range of $A_1$. Especially, we define a set $\mathbb{M}_1(2n,2k)$ by
\begin{equation}\label{A0}
\mathbb{M}_1(2n,2k): = Sp({2k},\mathbb{R}^{2n}) \cap \left\{ {\left. {\left. {\begin{bmatrix}
 \Phi & 0 \\
 0 & \Phi
 \end{bmatrix}} \right|\Phi  \in {\mathbb{R}^{n \times k}}} \right\}} \right ..
\end{equation}
If $A_1\in \mathbb{M}_1(2n,2k)$, $A_1 = {\rm{diag}}(\Phi, \Phi)$ for some $\Phi \in \mathbb{R}^{n\times k}$. Then, $A_1^TJ_{2n}A_1=J_{2k}$ if and only if $\Phi^T\Phi=I_k$, which implies that $\Phi$ is an element of the Stiefel manifold $V_k(\mathbb{R}^n)$. It follows that $\mathbb{M}_1(2n,2k)\cong V_k(\mathbb{R}^n)$, and
\begin{equation}
\mathbb{M}_1(2n,2k)= \left\{ {\left. {\left. {\begin{bmatrix}
 \Phi & 0 \\
 0 & \Phi
 \end{bmatrix}} \right|\Phi  \in {V_k(\mathbb{R}^n)}} \right\}} \right. .
\end{equation}

Let  ${R}$ and ${Q}$ denote two vector spaces; ${\rm{dim}}({R})=k$, ${\rm{dim}}({Q})=n$, and $k\le n$. Let $\mathbb{W}=T^*R$ and  $\mathbb{V}=T^*Q$. Suppose
 $f: {R}\to {Q}$ and  $\pi: {Q}\to {R}$ are linear mappings and satisfy $\pi\circ f={\rm{id}}_{ {R}}$.
  Let ${f_*}: {TR}\to {TQ}$ denote  the tangent lift of $f$. Let ${\left\langle   \cdot, \cdot \right\rangle}_q:=T_q^*Q\times T_qQ \to \mathbb{R}$ denote the natural pairing between tangent and cotangent vectors   at $q\in {Q}$, and let ${\left\langle  \cdot, \cdot  \right\rangle}_r:=T_r^*R\times T_rR \to \mathbb{R}$  denote the natural pairing between tangent and cotangent vectors   at $r\in {R}$.
The \emph{cotangent lift}, $f^*: \mathbb{V} \to \mathbb{W}$, of $f$ is a linear mapping that satisfies
\begin{equation}\label{colift}
{\left\langle {{f^*|_{q}}(p),w} \right\rangle _r} = {\left\langle {p,{f_{*r}}(w)} \right\rangle _q},
\end{equation}
     where $q=f(r)$, $w\in T_r{R}$, and $p\in T_q^*Q$.

  Using canonical coordinates, we have $r, w\in \mathbb{R}^k$, $q, p\in \mathbb{R}^n$. The linear mappings $f$ and $\pi$ are respectively denoted by matrices $\Phi \in \mathbb{R}^{n\times k}$ and $\Psi^T \in \mathbb{R}^{k\times n}$. Moreover, $\pi\circ f={\rm{id}}_R$ requires that $\Psi^T \Phi=I_k$. Since $f$ is a linear mapping, the tangent lift of $f$ at $r$ is represented by $f_{*r}(w)=\Phi w$.
  For $z\in \mathbb{W}$ and $x\in \mathbb{V}$, we have $z=[r;s]\in \mathbb{R}^{2k}$ and $x=[q;p]\in \mathbb{R}^{2n}$.
      Thus, (\ref{colift}) yields $s=\Phi^Tp$, and the cotangent lift $f^*: \mathbb{R}^{2n} \to \mathbb{R}^{2k}$ can be written as
  \begin{equation*}
  z=f^*(x)=B^+ x,
  \end{equation*}
   where $B^+ := {\rm{diag}}(\Psi^T, \Phi^T)$. Using (\ref{A++}),  we have $B=(B^+)^+=J_{2n}^TB^+J_{2k}={\rm{diag}} (\Phi,\Psi)$. Moreover, the constraint $\Psi^T\Phi=I_k$ implies that $B^+B=I_{2k}$. By Lemma \ref{lem:sym2}, we have $B\in Sp(2k, \mathbb{R}^{2n})$.

  Especially, when $\Phi=\Psi$,  $B$ degenerates to $A_1$, and the constraint $\Psi^T\Phi=I_k$ becomes $\Phi^T\Phi=I_k$. In this scenario, the range of $\Phi$ should approximately fit for both $q(t)$ and  $p(t)$. As Algorithm \ref{alg:lift} indicates, $\Phi$ can be computed by the SVD of an extended snapshot matrix $M_1\in \mathbb{R}^{n\times 2N}$, which is defined by
\begin{equation}\label{exsnap}
M_1:=[q(t_1), \ldots,  q(t_N), \gamma p(t_1), \ldots, \gamma p(t_N) ],
\end{equation}
where $\gamma$ represents a weighting coefficient. Let $\hat q(t)$ and $\hat p(t)$ denote approximating solutions based on a reduced system. We can choose $\gamma=\delta t$ if the goal is to minimize $\|\hat q(t)-q(t)\|_2$, and choose $\gamma=1$ if the goal is to minimize $\|[\hat q(t); \hat p(t)]-[q(t); p(t)]\|_2$.
In a similar way to the POD (SVD), the projection error of $M_x$ can be determined by the truncated singular values of $M_1$.

\begin{algorithm}
\caption{Cotangent Lift} \label{alg:lift}
\begin{algorithmic}
 \REQUIRE
An empirical data ensemble $\{q(t_i), p(t_i)\}_{i=1}^N$.
\ENSURE A symplectic matrix $A_1$ in block-diagonal form.
\STATE 1: Construct an extended snapshot matrix $M_1$ as (\ref{exsnap}).
\STATE 2: Compute the SVD of $M_1$ to obtain a POD basis matrix $\Phi$.
\STATE 3: Construct the symplectic matrix $A_1={\rm{diag}}(\Phi, \Phi)$.
\end{algorithmic}
\end{algorithm}

\medskip
\begin{theorem}
Suppose $M_x\in \mathbb{R}^{2n\times N}$ is the snapshot matrix defined by (\ref{snap}).  If we select $\gamma=1$ in (\ref{exsnap}), the symplectic matrix $A_1$ constructed by Algorithm \ref{alg:lift} is the optimal solution in $\mathbb{M}_1(2n,2k)$ that minimizes the error in the projection of $M_x$ onto the column space.
\end{theorem}
\medskip
\begin{proof}
Similar to (\ref{opt}), we can express the optimization problem as:
\begin{equation}\label{opt1}
\begin{aligned}
&{\rm{minimize}} \quad \ \|M_x-A_1A_1^+ M_x\|_F
\\ &{\rm subject \ to} \quad   A_1 \in \mathbb{M}_1(2n,2k).
\end{aligned}
\end{equation}
Let $M_q:=[q(t_1),\ldots,q(t_N)]$ and $M_p:=[p(t_1),\ldots,p(t_N)]$.  By definition, $M_x=[M_q;M_p]$. Moreover, $\gamma=1$ implies that  $M_1=[M_q,M_p]$. Since $A_1 \in \mathbb{M}_1(2n,2k)$, we have $A_1={\rm{diag}}(\Phi,\Phi)$ with ${\Phi ^T}\Phi = {I_k}$. Then, the objective function becomes
\begin{align*}
\|M_x-A_1A_1^+ M_x\|_F &=
{\left\| {\begin{bmatrix}
 M_q \\
 M_p
 \end{bmatrix} - \begin{bmatrix}
 \Phi & 0 \\
 0 & \Phi
 \end{bmatrix}\begin{bmatrix}
 {\Phi ^T} & 0 \\
  0 & {\Phi ^T}
 \end{bmatrix}\begin{bmatrix}
 M_q \\
 M_p
 \end{bmatrix}} \right\|_F} = {\left\| {\begin{bmatrix}
 (I_n - \Phi {\Phi ^T})M_q \\
 (I_n - \Phi {\Phi ^T})M_p
 \end{bmatrix}} \right\|_F}  \\
 & \! \! \! \! \! \! \! \! \! \! \! \! \! \! \! \! \! \! \! \! \! \! \! \! \! \! \! \! \!   \! \! \! \! \! \! \!\!  ={\left\| [(I_n - \Phi {\Phi ^T})M_q, \ (I_n - \Phi {\Phi ^T}) M_p] \right\|_F}={\left\| {(I_n - \Phi {\Phi ^T}){{\left[ M_q, \ M_p \right]}}} \right\|_F}={\left\| {M_1 - \Phi {\Phi ^T}M_1} \right\|_F}.
 \end{align*}
Thus, $\Phi$ can be directly solved by the truncated SVD of $M_1$,
\begin{equation}
M_1 \approx \Phi \Sigma \Psi^T,
\end{equation}
where the matrix $\Sigma$ is a $k\times k$ diagonal matrix with nonnegative real numbers on the diagonal; $\Phi$ and $\Psi$ are real matrices and satisfy $\Phi^T\Phi=\Psi^T\Psi=I_k$. Thus, the symplectic matrix $A_1$ constructed by Algorithm \ref{alg:lift} is the optimal solution for the optimization problem (\ref{opt1}).
 \hfill
\end{proof}
\medskip

It should be mentioned that in~\cite{Marsden:03h}, a tangent lift method is used to construct a reduced Euler-Lagrange equation to preserve the Lagrangian structure of the original system. Specifically, a POD basis matrix $\Phi\in \mathbb{R}^{n\times k}$ can be constructed by the SVD of a snapshot matrix $[q(t_1), \ldots, q(t_N)]$ for  $q(t)\in Q  \cong \mathbb{R}^n$.  Then, the original Lagrangian  $L(q,\dot q)$ in the tangent bundle ${TQ}$ is approximated by $\tilde L(r, \dot r)=L(\Phi r, \Phi \dot r)$ in ${TR}$, where $r(t)\in R  \cong \mathbb{R}^k$. Thus, a reduced system for $(r, \dot r)$ can be given by the Euler-Lagrange equation of $\tilde L(r, \dot r)$.

By the Legendre transformation, the reduced Lagrangian system can be transformed into a reduced Hamiltonian system. Meanwhile, the cotangent lift can yield another reduced Hamiltonian system. However, the two reduced systems are not equal in general, in two aspects.

 First, the two reduced systems lie on different subspaces. In either case, the subspace can be presented as the column space of $A_1={\rm{diag}}(\Phi,\Phi)$, where $\Phi$ is the POD basis matrix for the generalized coordinates.  In~\cite{Marsden:03h}, the tangent lift constructs $\Phi$ from a snapshot ensemble of $q(t)$. In this article, the proposed cotangent lift constructs $\Phi$ from a snapshot ensemble of $q(t)$ and $p(t)$, as (\ref{exsnap}) indicates.

 Second, the two reduced systems give different generalized momentums. Consider a  Lagrangian of the form $L(q, \dot q)={1\over 2} \dot q^TM\dot q-V(q)$, where $M\in \mathbb{R}^{n\times n}$ denotes  mass matrix and $V:\mathbb{R}^n \to \mathbb{R}$ denotes the potential function.  Then, the reduced Lagrangian is given by $\tilde L(r, \dot r)={1 \over 2} \dot r^T \tilde M  \dot r-V(\Phi r)$, where $\tilde M=\Phi^T M \Phi\in \mathbb{R}^{k\times k}$ denotes the reduced mass matrix. The Legendre transform produces the reduced momentum by $s_1= \tilde M \dot r=\Phi^T M \Phi \Phi^T \dot q$; while the cotangent lift approach produces the reduced momentum by $s_2=\Phi^T p=\Phi^T M \dot q$.  Unless $\dot q(t)$ resides on  $\rm{Range}(\Phi)$ or $M=I_n$, $s_1= s_2$  does not hold  in general.

\subsection{Complex SVD}\label{sec:csvd}
This section proposes an SVD-based algorithm to construct a symplectic basis matrix,  such that the off-diagonal blocks are non-zero submatrices. If we use $q(t)+\iota p(t)$ to describe the solution trajectory in the phase space, we can construct a complex snapshot matrix  $M_2\in \mathbb{C}^{n\times N}$ by
\begin{equation}\label{Z}
M_2:=[q(t_1)+ \iota p(t_1), \ldots, q(t_N), +\iota p(t_N) ].
\end{equation}
By definition, we have $M_2=M_q+ \iota M_p.$  Suppose a unitary matrix $U\in \mathbb{C}^{n\times k} $ minimizes the error in the projection of $M_2$ onto the column space. Then, $U$ is the solution of the following optimization problem:
\begin{equation}
\begin{aligned}
 & {\rm{minimize }}  \quad \ {\left\| {M_2 - U {U ^H}M_2} \right\|_F} \\
 & {\rm{subject \ to  }} \quad {U ^H}U = {I_k}.\\
 \end{aligned}
 \end{equation}
Here $U^H$ is the conjugate transpose of $U$. In fact, $U$ can be obtained by the truncated SVD of $M_2$,
\begin{equation}\label{svd}
M_2 \approx U\Sigma V^H,
\end{equation}
where the matrix $\Sigma$ is a $k\times k$ diagonal matrix with nonnegative real numbers on the diagonal, and $U$ and $V$ are complex matrices and satisfy $U^HU=V^HV=I_k$.  Let $V_k(\mathbb{C}^n)$ denote the Stiefel manifold in $\mathbb{C}^n$. Then, its element $U\in V_k(\mathbb{C}^n)$ has the form  $U=\Phi+\iota \Psi$, where $\Phi, \Psi \in \mathbb{R}^{n \times k}$. We define a mapping $\mathscr{A}: V_k(\mathbb{C}^n) \to \mathbb{R}^{2n\times 2k}$ by the formula
\begin{equation}\label{iotau}
\mathscr{A} (U)=  \begin{bmatrix}
 {\Phi} & -{\Psi} \\
 {\Psi}& {\Phi}
 \end{bmatrix}.
\end{equation}

\begin{lemma}\label{lem:bij}
The mapping $\mathscr{A}$ is injective. The image of $\mathscr{A}$ is equal  to  $\mathbb{M}_2(2n,2k)$, where
 \begin{equation}
\mathbb{M}_2(2n,2k):=Sp({2k},\mathbb{R}^{2n})  \cap \left\{ {\left. {\left. {\begin{bmatrix}
 {\Phi} & -{\Psi} \\
 {\Psi} & {\Phi}
 \end{bmatrix}} \right|\Phi,\Psi  \in {\mathbb{R}^{n\times k}}} \right\}} \right. .
 \end{equation}
\end{lemma}
\medskip

\begin{proof} It follows from $\mathscr{A}$'s definition that it is injective.  If $\Phi+\iota \Psi\in {V_k(\mathbb{C}^n)}$, then $(\Phi+\iota \Psi)^H(\Phi+\iota \Psi)=I_k$, which is equivalent to
\begin{equation}\label{comcon}
\Phi^T\Phi+\Psi^T\Psi=I_k, \ \ \ \Phi^T\Psi=\Psi^T\Phi.
\end{equation}
Let $A_2=\mathscr{A} (\Phi+\iota \Psi)$. Using (\ref{comcon}), it is easy to verify that $A_2^TJ_{2n}A_2=J_{2k}$. Thus, $A_2 \in \mathbb{M}_2(2n,2k)$, i.e., $\mathscr{A}({V_k(\mathbb{C}^n)}) \subset \mathbb{M}_2(2n,2k)$.

Conversely, if $A_2\in \mathbb{M}_2(2n,2k)$, then $
A_2^TJ_{2n}A_2=J_{2k}$.
Moreover, we can write $A_2=[\Phi, -\Psi; \Psi, \Phi]$ for some $\Phi,\Psi\in \mathbb{R}^{n\times k}$. Plugging it into $A_2^TJ_{2n}A_2=J_{2k}$ gives (\ref{comcon}). It follows that, $(\Phi+\iota \Psi)^H(\Phi+\iota \Psi)=I_k$. As a result, $\Phi+\iota \Psi\in {V_k(\mathbb{C}^n)}$, and $\mathscr{A} ^{-1}(\mathbb{M}_2(2n,2k))\subset {V_k(\mathbb{C}^n)}$.
\hfill
\end{proof}

\begin{algorithm}
\caption{Complex SVD} \label{alg:complex}
\begin{algorithmic}
 \REQUIRE
An empirical data ensemble $\{q(t_i), p(t_i)\}_{i=1}^N$.
\ENSURE A symplectic matrix $A_2$ in block form.
\STATE 1: Construct a complex snapshot matrix $M_2$ as (\ref{Z}).
\STATE 2: Compute the SVD of $M_2$ to obtain a basis matrix $\Phi+\iota \Psi$.
\STATE 3: Construct the symplectic matrix $A_2=[\Phi, -\Psi; \Psi, \Phi]$.
\end{algorithmic}
\end{algorithm}

Lemma \ref{lem:bij} implies that a symplectic matrix $A_2$ can be constructed through the mapping $\mathscr{A}$. Algorithm \ref{alg:complex} outlines the procedure. Since both $\mathscr{A}$ and $\mathscr{A}^{-1}$ are smooth mappings,  $\mathbb{M}_2(2n,2k)\cong V_k(\mathbb{C}^n)$, and $\mathbb{M}_2(2n,2k)$ is a submanifold of $Sp({2k},\mathbb{R}^{2n})$.

Moreover, by substituting $A_2 = [\Phi, \Psi; -\Psi, \Phi]$ into (\ref{A+}), we obtain $A_2^+=J_{2k}^TA^T_2J_{2n}= A_2^T$. It follows that $A_2^TA_2=A_2^+A_2=I_{2k}$, i.e., $A_2\in V_{2k}(\mathbb{R}^{2n})$. Conversely, for any $A_2=[\Phi, \Psi; -\Psi, \Phi] \in V_{2k}(\mathbb{R}^{2n})$, (\ref{comcon}) holds, which means that $A_2 \in Sp({2k},\mathbb{R}^{2n})$.  Therefore, an equivalent definition of  $\mathbb{M}_2(2n,2k)$ is given by
    \begin{equation}\label{b514}
    \mathbb{M}_2(2n,2k)=V_{2k}(\mathbb{R}^{2n})  \cap \left\{ {\left. {\left. {\begin{bmatrix}
 {\Phi} & -{\Psi} \\
 {\Psi} & {\Phi}
 \end{bmatrix}} \right|\Phi,\Psi  \in {\mathbb{R}^{n\times k}}} \right\}} \right. .
 \end{equation}

\medskip
\begin{lemma}\label{lem:iso}
$V_k(\mathbb{C}^n)$ is isomorphic to $Sp(2k, \mathbb{R}^{2n})\cap  V_{2k}(\mathbb{R}^{2n})$.
\end{lemma}
\medskip

\begin{proof}
Let $A_q=[\xi_1,\ldots, \xi_k]$, $A_p=[\zeta_1, \ldots, \zeta_k]$, and $A_2=[A_q, A_p]\in Sp(2k, \mathbb{R}^{2n})\cap  V_{2k}(\mathbb{R}^{2n})$.  For any $i\in \{1, \ldots, k\}$, we have $\|\xi_i\|=\|\zeta_i\|=1$, and $\Omega(\xi_i, \zeta_i)=1$. It follows that $\|J_{2n}\xi_i\|=1$, and the inner product $\left\langle J_{2n}\xi_i, \zeta_i \right\rangle=1$. The Cauchy--Schwarz inequality states that $\left\langle J_{2n}\xi_i, \zeta_i \right\rangle\le \|J_{2n}\xi_i\| \cdot \|\zeta_i\|$, and  two sides are equal if and only if $J_{2n}\xi_i$ and $\zeta_i$ are parallel. Thus, we must have $J_{2n}\xi_i=\zeta_i$.  It follows that $A_2$ must have the block form $[A_q, J_{2n}A_q]$, or $[\Phi, -\Psi; \Psi, \Phi]$ if $A_q$ is written as $[\Phi; \Psi]$. Therefore, $A_2 \in \mathbb{M}_2(2n, 2k)$. By the proof in Lemma \ref{lem:bij}, we have $\Phi+\iota \Psi \in V_{k}(\mathbb{C}^n)$.

Conversely, if $\Phi+\iota \Psi \in V_{k}(\mathbb{C}^n)$, the mapping (\ref{iotau}) yields  $\mathscr{A}(\Phi+\iota \Psi) \in Sp(2k, \mathbb{R}^{2n})$ and $\mathscr{A}(\Phi+\iota \Psi) \in V_{2k}(\mathbb{R}^{2n})$.
 \hfill
\end{proof}

\medskip

Notice that  $\mathscr{A}$ also  preserves algebraic structures, as one can easily verify the following lemma.

\medskip
\begin{lemma}\label{lem:CD}
Let $C\in \mathbb{C}^{n_1\times n_2}$ and $D\in \mathbb{C}^{n_2\times n_3}$. Then, we have
$\mathscr{A}(C)\mathscr{A}(D)=\mathscr{A}(CD)$ and $\mathscr{A}(C^H)=(\mathscr{A}(C))^T$.
\end{lemma}
\medskip

\begin{theorem}
Suppose $M_x\in \mathbb{R}^{2n\times N}$ is the snapshot matrix defined by (\ref{snap}).  The symplectic matrix $A_2$ constructed by Algorithm \ref{alg:complex} is the optimal solution in $\mathbb{M}_2(2n,2k)$ that minimizes the error in the projection of $[M_x, -J_{2n}M_x]$ onto the column space.
\end{theorem}
\medskip

\begin{proof}
By Lemma \ref{lem:CD}, the truncated SVD of $M_2$ given by (\ref{svd}) yields
\begin{equation}\label{svd1}
\mathscr{A}(M_2) \approx \mathscr{A}(U) \mathscr{A}(\Sigma) \mathscr{A}(V^H)= \mathscr{A}(U) \mathscr{A}(\Sigma) (\mathscr{A}(V))^T.
\end{equation}
Since $U^HU=I_k$, by Lemma \ref{lem:CD}, we have
\begin{equation*} (\mathscr{A}(U))^T\mathscr{A}(U)=\mathscr{A}(U^H)\mathscr{A}(U)=\mathscr{A}(I_k)=I_{2k}.
\end{equation*}
Similarly, $(\mathscr{A}(V))^T\mathscr{A}(V)=I_{2k}$ holds due to $V^HV=I_k$. Moreover, $\mathscr{A}(\Sigma)$ is a real diagonal matrix that contains the first $2k$ dominant singular values of $\mathscr{A}(M_2)$. Thus, (\ref{svd1}) provides the truncated SVD for $\mathscr{A}(M_2)$.

In Algorithm \ref{alg:complex}, the symplectic matrix is constructed   by $A_2=\mathscr{A}(U)$. Meanwhile, using the definition of $M_2$ and $M_x$, we have
$$\mathscr{A}(M_2)=
\begin{bmatrix}
 M_q &  -M_p \\
 M_p & M_q \\
\end{bmatrix}
 =\begin{bmatrix}
 M_x, \ -J_{2n}M_x
 \end{bmatrix}.
 $$
Therefore, $A_2$ minimizes the truncation error due to the construction of $[M_x, -J_{2n}M_x]$ using a symplectic subspace with a fixed dimension.
  \hfill
\end{proof}
\medskip

It should be emphasized that because the complex SVD is designed to fit $[M_x, -J_{2n}M_x]$, rather than $M_x$ itself. As a result, Algorithm \ref{alg:complex} can only construct a near optimal matrix in $\mathbb{M}_2$.

\subsection{Nonlinear Programming (NLP)} Although it is often too expensive to solve the optimization problem (\ref{opt}) directly, one can search a near optimal solution over a subset of $Sp({2k},\mathbb{R}^{2n})$. Especially, if one has a pre-specified basis matrix $A_1 \in \mathbb{M}_1(2n,2r)$ or $A_2 \in \mathbb{M}_2(2n,2r)$, with $k\le r\le n$,  one may assume that the near optimal solution $A_3\in Sp({2k},\mathbb{R}^{2n})$ is a linear transformation of $A_1$ or $A_2$.

 Now suppose the cotangent lift yields a symplectic matrix $A_1$ in $\mathbb{M}_1(2n, 2r)$. If  ${\rm{Range}}(A_3) \subset {\rm{Range}}(A_1)$, we have
 \begin{equation}\label{AC}
 A_3=A_1 \cdot C,
 \end{equation}
 where $C\in \mathbb{R}^{2r\times 2k}$ is the coefficient matrix of $A_3$ with respect to the basis vectors of $A_1$.   Plugging (\ref{AC}) into $A_3^TJ_{2n}A_3=J_{2k}$ and using $A_1^TJ_{2n}A_1=J_{2r}$ give
 \begin{equation}\label{CJC}
 C^T J_{2r} C=J_{2k},
 \end{equation}
 which implies $C \in Sp({2k},\mathbb{R}^{2r})$. Notice $(A_1C)^+=C^+A_1^+$, the original optimization problem (\ref{opt}) reduces to
 \begin{equation}\label{opt2}
\begin{aligned}
&{\rm{minimize}} \quad \ \|M_x-A_1CC^+A_1^+ M_x\|_F \\
 &{\rm subject \ to}  \quad  C^TJ_{2r}C=J_{2k}.
\end{aligned}
\end{equation}
Let $A_1=\rm{diag}(\Phi, \Phi)$, where $\Phi\in \mathbb{R}^{n\times r}$. The initial value for (\ref{opt}) could be $A_3=\rm{diag}(\Phi', \Phi')$, where $\Phi'$ denotes the first $k$ columns of $\Phi$. Correspondingly, the initial value for (\ref{opt2}) is given by $C={\rm{diag}}(I_{r\times k}, I_{r\times k})$, where $I_{r\times k}$ denotes the first $k$ columns of the identity matrix $I_r$. Since (\ref{opt2}) requires the optimization of the coefficient matrix $C$ over a smaller domain to fit the empirical data, the computational cost is significantly lower than the original optimization problem (\ref{opt}) when $r\ll n$. Finally, it should be mentioned that the proposed NLP algorithm here is analogous to the optimization algorithm in~\cite{WillcoxK:10a}, where an optimal POD basis matrix is constructed from a linear transformation of the snapshot matrix.

 \begin{algorithm}
\caption{Nonlinear Programming} \label{alg:nlp2}
\begin{algorithmic}
 \REQUIRE
An empirical data ensemble $\{q(t_i), p(t_i)\}_{i=1}^N$.
\ENSURE A symplectic matrix $A_3 \in Sp({2k},\mathbb{R}^{2n})$.
\STATE 1: Construct a symplectic matrix $A_1 \in \mathbb{M}_1(2n,2r)$ with $r>k$ by the cotangent lift.
\STATE 2: Solve (\ref{opt2})   and obtain a coefficient matrix $C \in Sp({2k},\mathbb{R}^{2r})$.
\STATE 3: Construct the symplectic matrix $A_3=A_1 \cdot C$.
\end{algorithmic}
\end{algorithm}

So far, three different algorithms have been proposed to construct a symplectic basis matrix. Corresponding to three manifolds with the inclusion maps:
$$V_k(\mathbb{R}^n)  \xhookrightarrow{} V_k(\mathbb{C}^n) \xhookrightarrow{\mathscr{A}} Sp(2k, \mathbb{R}^{2n}),$$
we propose the cotangent lift,  complex SVD, and  NLP. The cotangent lift and complex SVD algorithms are faster and more easily implemented in offline computation; their computational costs only involve the SVD. However, both algorithms search optimal basis matrices in submanifolds of $Sp(2k, \mathbb{R}^{2n})$, rather than in $Sp(2k, \mathbb{R}^{2n})$ itself. Therefore, they sacrifice certain accuracy to fit the empirical data in order to reduce costs. By contrast, the NLP is more expensive in offline computation, since it requires solving an optimization problem in $Sp(2k, \mathbb{R}^{2r})$ based on a pre-specified basis matrix constructed by another algorithm. However, the NLP can result in a symplectic matrix to fit the empirical data with less projection error.

For the cotangent lift, we have a parameter $\gamma$ in (\ref{exsnap}) to balance the truncation error due to the construction of $p(t)$ and $q(t)$ using a symplectic subspace. At first glance, the other two algorithms do not have a similar weighting option. However, we can always construct  a linear transformation from $x=[q; p]$  to $x_\gamma=[q; p_\gamma]$ by $p_\gamma=\gamma p$, and then solve  the rescaled Hamiltonian equation based on $\tilde H(q, p_\gamma): =H(q, p_\gamma/\gamma)$. The fully-resolved  rescaled Hamiltonian system is equivalent to the original one; depending on the subspace on which the reduced system lives, however, the reduced models for the original and rescaled systems are not equivalent in general. A weighted data ensemble for the rescaled system can be defined as
\begin{equation}\label{gamma}
M_{x_\gamma}:=[x_\gamma(t_1), \ldots, x_\gamma(t_N)].
 \end{equation}
 Then, a symplectic subspace can be constructed to fit $M_{x_\gamma}$ by any of the aforementioned PSD algorithms. Thus, the complex SVD and NLP can also flexibly balance the truncation error of  $p(t)$ and $q(t)$ by choosing a suitable value of $\gamma$.

\section{Symplectic Model Reduction of Nonlinear Hamiltonian Systems}\label{sec:deim}
As an approximation of the symplectic Galerkin projection, the SDEIM is developed in this section. The motivation of the SDEIM is to reduce the computational complexity of a nonlinear Hamiltonian system while simultaneously preserving the symplectic structure. Before introducing the SDEIM, we will give a review of the classical DEIM.

\subsection{Discrete Empirical Interpolation Method (DEIM)}\label{sec:DEIM}
Let $x\in \mathbb{R}^n$ denote the state variable in the original space and let $f: \mathbb{R}^n \to \mathbb{R}^n$ denote the discretized vector field.
The full-order dynamical system can be described by an initial value problem
\begin{equation} \label{general1}
\dot x = f(x)=L x+ f_N(x); \quad x(0)=x_0,
\end{equation}
where the original vector field $f(x)$ is split into a linear part $L x$ with $L\in \mathbb{R}^{n\times n}$ and a nonlinear part $f_N(x)$ with $f_N: \mathbb{R}^n \to \mathbb{R}^n$. \footnote{In general, the choice of $L x$ and $f_N(x)$ is not unique. If we let $L=0$, then $f_N(x)=f(x)$. In this scenario, one can avoid computing $\tilde L z$  and save some computational cost in the online stage. However, since the DEIM is only an approximation for the standard Glakin projection, it inevitably introduces extra error to evaluate the linear term when $L x$ is absorbed in $f_N(x)$. Therefore, it is desired to separate the original vector filed $f(x)$ and compute the linear term $L x$ via Galekin projection, especially when the  $Lx$ takes the dominance over $f_N(x)$.}

Let $\Phi \in \mathbb{R}^{n\times k}$ denote a POD basis matrix. Then, the Galerkin projection can be used to obtain a reduced system on the column space of $\Phi$,
\begin{equation} \label{rom1}
\dot z  =   \Phi^T f(\Phi z)= \tilde L z + \Phi^T f_N (\Phi z); \quad    z_0=\Phi^T x_0,
\end{equation}
where  $z(t)\in \mathbb{R}^k$ is the reduced state, and $\tilde L=\Phi^T L \Phi\in \mathbb{R}^{k\times k}$ is the reduced linear operator.

 According to some previous studies~\cite{SorensenDC:10a,PetzoldLR:03a}, the POD-Galerkin can achieve computational savings only when the analytical formula of the nonlinear vector term $\Phi^T f_N(\Phi z)$ can be simplified, especially if $f_N(x)$ is a low-order polynomial in $x$. Otherwise, one usually needs to compute the state variable $\hat x:=\Phi z$ in the original coordinate system, evaluate the nonlinear vector field $f_N(\hat x)$, and then project $f_N(\hat x)$ onto the column space of $\Phi$. In this scenario, solving the POD reduced system (\ref{rom1}) could be more expensive than solving the original full-order system (\ref{general1}).

 DEIM focuses on approximating $f_N$ so that a certain coefficient matrix can be pre-computed and, as a result, the complexity in evaluating $f_N$ becomes proportional to the small number of selected spatial indices~\cite{SorensenDC:10a}. Let $\beta = [\beta_1; \ldots; \beta_{m}] \in \mathbb{R}^{m}$ be an index vector, and $\beta_i \in \{1, \ldots, n\}$ a index. Define an ${n\times m}$ matrix
\begin{equation}\label{proj0}
P := [e_{\beta_1}, \ldots, e_{\beta_m}],
\end{equation}
where $e_{\beta_i}$ is the $\beta_i$th column of the identity matrix $I_{2n}$.   Then, left multiplication of $f_N(x)$ with $P^T$ projects $f_N(x)$ onto $m$ elements corresponding to the index vector $\beta$.
Now suppose $f_N(x)$ resides approximately on the range of an $n\times m$ matrix $\Psi$, then there exists a corresponding coefficient vector $\tau\in \mathbb{R}^{m}$ such that $f_N(x) \approx \Psi \tau(x)$. The coefficient vector $\tau(x)$ can be determined by matching the $f_N(x)$ at selected $m$ spatial indices, i.e., $P^T  f_N(x) = P^T \Psi\tau(x)$. Suppose $P^T f_N(x)$ is nonsingular. Then, we have $\tau(x)=(P^T \Psi)^{-1}P^T f_N(x)$. Thus, the approximation $\hat f_N(x)$ of the nonlinear vector term $f_N(x)$ becomes
\begin{equation}\label{approxnon1}
 \hat f_N(x) = \Psi \tau(x) = \Psi (P^T \Psi)^{-1}P^T  f_N(x),
\end{equation}
and the reduced system (\ref{rom1}) can be approximated as
\begin{equation} \label{romdeimf}
\dot z =   \tilde L z +W  g(z),
\end{equation}
where $W=\Phi^T \Psi (P^T \Psi)^{-1}$, and $g(z)=P^T  f_N(\Phi z)$. Notice that $W$ is calculated only once at the offline stage. At the online stage, $g(z)$ is evaluated on $m$ spatial indices of $f_N(\Phi z)$.  Therefore, the complexity of the DEIM-reduced system (\ref{romdeimf}) could be independent of dimension $n$ of the original system.

\medskip
\begin{algorithm}
\caption{Greedy algorithm to construct an index vector $\beta$} \label{alg:greedy}
\begin{algorithmic}
\REQUIRE A basis matrix $\Psi=[\psi_1, \ldots, \psi_{m}]\in \mathbb{R}^{n\times m}$.
\ENSURE An index vector $\beta=[\beta_1; \ldots ;\beta_{m}]\in \mathbb{R}^{m}$.
\STATE 1: Select the first interpolation index
$[ \rho ,\beta_1]=\max\{|\psi_1|\}$.
\STATE 2: Initialize $U=[\psi_1]$, $\beta=\beta_1$.
 \FOR{$i = 2$ to $m$}
 \STATE 3: Solve the coefficient vector $\tau$ to match $\psi_i$,
 $U(\beta,:)\tau= \psi_i(\beta)$.
 \STATE 4: Calculate the residual $r=\psi_i-U \tau$.
 \STATE 5: Select the interpolation index corresponding to the largest magnitude of the residual $r$,
 $[\rho ,\beta_i]=\max\{|r|\}$.
 \STATE 6: Update $U = [U, \ \psi_i]$, $\beta =[\beta; \beta_i]$. \ENDFOR
\end{algorithmic}
\end{algorithm}

In order to construct (\ref{romdeimf}), the SVD can be applied to construct the POD basis matrix $\Phi$ based on an empirical data ensemble $[x(t_1), \ldots, x(t_N)]$ and the collateral POD basis matrix $\Psi$ based on another data ensemble $[f_N(x(t_1)), \ldots, f_N(x(t_N))]$ for the nonlinear vector term.
Moreover, a greedy algorithm can be applied to construct the index vector $\beta$~\cite{SorensenDC:10a}, as listed in Algorithm \ref{alg:greedy}.\footnote{The MATLAB notations $B(\beta,:)$ and $a_i(\beta)$ are used here to represent the operation of selecting rows out of a matrix (or a vector).} Initially, we select the first interpolation index $\beta_1 \in \{1,\ldots,n\}$ corresponding to the first basis function $\psi_1$ with the largest magnitude. The remaining interpolation indices, $\beta_i$ for $i= 2, \ldots,m$, respectively correspond to the largest magnitude of the residual $r$, where $r$ is the residual between the input basis $\psi_i$ and its projection onto the column space of $U$. Especially, In Step 5, $[\rho ,\beta_i]=\max\{|r|\}$ means $\rho=|r(\beta_i)|={\rm{max}}_{j=1,\ldots,n}{|r(j)|}$. In Step 6, we add a column vector $\psi_i$ (and an element $\beta_i$) to a matrix $U$ (and a vector $\beta$). It has been proven that $\rho  \ne 0$ implies that $P^TA$ is nonsingular~\cite{SorensenDC:10a}. Thus, the DEIM approximation of the nonlinear vector term $\hat f_N(x)$ in (\ref{approxnon1}) is well-defined.

\subsection{Symplectic Discrete Empirical Interpolation Method (SDEIM)}
Similar to (\ref{general1}), the original Hamiltonian can also be split into a linear part and a nonlinear part, i.e., $H(x)=H_1(x)+H_2(x)$, such that $\nabla_x H_1(x)=Lx$ for a real symmetric matrix $L$, and $\nabla_x H_2(x)=f_N(x)$ for a nonlinear function $f_N$. Thus, the original Hamiltonian system is given by
\begin{equation}\label{split}
\dot x=J_{2n} \nabla_x H(x)=Kx+J_{2n}f_N(x),
\end{equation}
where $K =J_{2n}L\in \mathfrak{sp}(\mathbb{R}^{2n})$. Analogous to (\ref{rom1}), the symplectic Galerkin projection yields the following reduced Hamiltonian system
\begin{equation} \label{nonlinear}
\dot z = {A^ + }(Kx+J_{2n}f_N(x))=  \tilde K z +J_{2k}A^Tf_N(Az),
\end{equation}
where $A\in Sp(2k, \mathbb{R}^{2n})$.
Thus, unless $A^Tf_N(Az)$ can be analytically simplified, the computational complexity of (\ref{nonlinear}) still depends on $2n$. In order to save the computational cost, one can use the DEIM approximation (\ref{approxnon1}) to approximate the nonlinear vector term $f_N$.
Let $\Psi\in \mathbb{R}^{2n \times m}$ denote the collateral POD basis for $f_N(x)$, and let $P\in \mathbb{R}^{2n \times m}$ denote the projection matrix with the form (\ref{proj0}).
Then,    (\ref{nonlinear}) can be approximated as
\begin{equation} \label{romdeim}
\dot z =   \tilde K z +J_{2k} W g(z),
\end{equation}
where $W=A^T  \Psi (P^T \Psi)^{-1}$ and $g(z)=P^T f_N(Az)$. Strictly speaking, (\ref{romdeim}) is not necessarily to be Hamiltonian. However, when the DEIM offers a good approximation for the PSD reduced system (\ref{rom}), one may expect that evolving (\ref{romdeim}) will not yield a large energy variation.

\medskip
\begin{definition}
The {\textbf{SDEIM}} of a nonlinear Hamiltonian system $\dot x=J_{2n}\nabla_x H(x)$, or $\dot x =  K x +J_{2n} f_N(x)$, with an initial condition $x(0)=x_0$ is given by (\ref{romdeim}), with the initial condition $z_0=A^+x_0$.
\end{definition}
 \medskip

  Both the cotangent lift (in Section \ref{sec:clm}) and the complex SVD (in Section \ref{sec:csvd}) can be used to construct a symplectic matrix $A\in \mathbb{M}_2(2n, 2k)$ based on an empirical data ensemble. Notice that (\ref{b514}) implies that $\mathbb{M}_2(2n,2k)\subset V_{2k}(\mathbb{R}^{2n})$. Thus, if we choose $\Psi$ such that
\begin{equation} \label{suf}
A=\Psi\in \mathbb{M}_2(2n,2k),
\end{equation}
then $A^T \Psi=I_{2n}$. It follows that $W=(P^TA)^{-1}$.  Since $x(t)\in \mathbb{V}$ and  $\nabla_x  H_2(x)=f_N(x)$, $[x(t); f_N(x(t))]$ is a trajectory  in $T^*\mathbb{V}$. By assuming $A=\Psi$ in (\ref{suf}), we actually lift a mapping $\sigma:\mathbb{W}\to \mathbb{V}$ to $\sigma_*:T\mathbb{W} \to T\mathbb{V}$ via a $4n\times 4k$ matrix, ${\rm{diag}} (A,A)$. Using a similar idea from Section \ref{sec:clm}, $A$ can be constructed by an extended data ensemble,
\begin{equation}\label{exsnap1}
M_3:=[x(t_1), \ldots, x(t_N), f_N(x(t_1)), \ldots, f_N(x(t_N)) ],
\end{equation}
that contains both the state $x(t_i)$ and the nonlinear term $f_N(x(t_i))$.

Regarding the computational complexity of the SDEIM, $\tilde K$ and $W$ are calculated only once at the beginning. For each step in the online stage, the nonlinear vector term $g(z)$ is only evaluated on selected $2k$ spatial indices of $f_N(Az)$. Thus, the complexity of SDEIM is also $O(1)$  when $k$ and $m'$ are fixed. Here $m'$ denotes the number of elements of $Az$ that are required to compute the $2k$ spatial indices of $P^T f_N(Az)$.

Table \ref{tab:podvspsd} compares the POD-Galerkin with the proposed symplectic model reduction; it serves as a short summary of Sections \ref{sec:projection}-\ref{sec:deim}.

\begin{table}[H]
\begin{center}
\caption{The POD-Galerkin vs. the symplectic model reduction.}
 \label{tab:podvspsd}
\begin{tabular}{|l|l |l|}
\hline
&  POD-Galerkin approach  & Symplectic model reduction
\\
\hline
Original system & $\dot x=f(x)$ with $x\in \mathbb{R}^n$ & $\dot x=J_{2n}\nabla_x H(x)$ with  $x\in \mathbb{R}^{2n}$ \\
 \hline
 Reduced state & \begin{tabular}{@{}l@{}} Orthogonal projection: \\   $z=\Phi^T x\in \mathbb{R}^k$ \end{tabular}
    & \begin{tabular}{@{}l@{}} Symplectic projection: \\   $z=A^+ x\in \mathbb{R}^{2k}$ \end{tabular}  \\ \hline
 Reduced system & \begin{tabular}{@{}l@{}} Galerkin projection: \\   $\dot z =\Phi^T f(\Phi z)$ \end{tabular} &
             \begin{tabular}{@{}l@{}} Symplectic Galerkin projection: \\   $\dot z = J_{2k}{\nabla _z}H(Az)$\end{tabular}       \\
  \hline
 \begin{tabular}{@{}l@{}} Properties of \\ reduced system \end{tabular} &  No stability guarantee &  \begin{tabular}{@{}l@{}} Energy preservation  \\ Stability preservation \end{tabular}\\
 \hline
 Basis matrix & Orthonormal: $\Phi^T \Phi=I_k$ & Symplectic: $A^TJ_{2n}A=J_{2k}$\\
 \hline
 \begin{tabular}{@{}l@{}}
 Domain of \\ basis matrix
 \end{tabular} & \begin{tabular}{@{}l@{}} Stiefel manifold  $V_k(\mathbb{R}^n)$ \end{tabular}& \begin{tabular}{@{}l@{}} Symplectic Stiefel manifold \\ $Sp({2k},\mathbb{R}^{2n})$ \end{tabular}\\
 \hline
 \begin{tabular}{@{}l@{}} Constructing \\ basis  matrix  \end{tabular}
 & \begin{tabular}{@{}l@{}} Proper orthogonal \\ decomposition (POD) \end{tabular}& \begin{tabular}{@{}l@{}} Proper symplectic  decomposition:  \\ (PSD)\\(a) Cotangent lift  \\ (b) Complex SVD \\ (c) Nonlinear programming (NLP) \end{tabular}\\
 \hline
 \begin{tabular}{@{}l@{}} Simplifying \\ nonlinear  terms \end{tabular} &   DEIM: Equation (\ref{romdeimf})   & SDEIM: Equation (\ref{romdeim})
 \\
 \hline
\end{tabular}
\end{center}
\end{table}

\section{Numerical Examples} \label{sec:numerical}
 In this section, the performance of symplectic model reduction is illustrated in numerical simulation of wave equations. After deriving the Hamiltonian form of general wave equations, we first study a linear wave equation numerically and focus on demonstrating the capability of PSD algorithms to deliver stability-preserving reduced systems. Then we simulate the nonlinear sine-Gordon equation to illustrate that the SDEIM is able to deliver accurate and long-time stable results with significant speedups.

\subsection{Hamiltonian Formulation for Wave Equations}
Let $u=u(t,x)$.  Consider the one-dimensional semi-linear wave equation with constant moving speed $c$ and a nonlinear vector term $g(u)$,
\begin{equation}\label{lwave}
  u_{tt}=  c^2   u_{xx}-g(u),
\end{equation}
on space $x\in[0,l]$. With the generalized coordinates $q=u$ and the generalized momenta $p= u_t$, the Hamiltonian PDE associated with (\ref{lwave}) is given by
\begin{equation}\label{qp}
 \dot q= \frac  {\delta {H}}{\delta p}, \qquad
  \dot p= -\frac {\delta  {H}}{\delta q},
\end{equation}
where the Hamiltonian is defined as
  \begin{equation}\label{hamiqp}
 H(q,p) = \int_0^l {dx\left[ {\frac{1}{2}{p^2}{\rm{ + }}\frac{{\rm{1}}}{{\rm{2}}}{c^2 q_x^2}} +G(q) \right]}, \qquad G'(q)=g(q).
  \end{equation}

A fully resolved model of (\ref{qp}) can be constructed by a structure-preserving finite difference discretization~\cite{BridgesTJ:06a}. With $n$ equally spaced grid points, the spatial discretized Hamiltonian is given by
\begin{equation}\label{waveeng}
{H_d}(y) = \sum\limits_{i = 1}^n {\Delta x\left[ \frac{1}{2}p_i^2 + \frac{c^2{({q_{i + 1}} - {q_i})}^2}{4\Delta {x^2}}+ \frac{c^2{({q_{i}} - {q_{i-1}})}^2}{4\Delta {x^2}} +G(q_i) \right]},
\end{equation}
where $q_i := u(t,x_i)$, $p_i := u_t(t,x_i)$,  $y: =[q_1; \ldots; q_n; p_1; \ldots; p_n]$, and $x_i = i\Delta x$. In the limit $\Delta x \to 0$ and $n\Delta x=l$, (\ref{waveeng}) converges to (\ref{hamiqp}). Now, the full model is represented by a Hamiltonian ODE system,
\begin{equation}\label{wavedis}
\frac{{{\rm{d}}{y}}}{{{\rm{d}}t}} = {J_d}\nabla_y {H_d}, \qquad {J_d} = \frac {J_{2n}}{\Delta x}.
\end{equation}
Let $D_{xx}\in \mathbb{R}^{n\times n}$ denote the  the three-point central difference approximation for the spatial derivative $\partial_{xx}$. We define a Hamiltonian matrix by
 \begin{equation}\label{wavematrix}
 K=
 \begin{bmatrix}
 0_n & I_n \\  c^2 {D}_{xx} & 0_n
 \end{bmatrix}.
 \end{equation}
 Then, (\ref{wavedis}) can be written in the form
\begin{equation}\label{expansion}
\dot y=Ky+J_{2n}f_N(y),
\end{equation}
 where the nonlinear vector term $f_N(y)$ is a vector in $\mathbb{R}^{2n}$ with zeros in the last $n$ elements. We have $f_N(y)=[g(q);0_{n\times 1}]$ for periodic and Neumann boundary conditions, and $f_N(y)=[g(q)- {{c^2} \over {\Delta x^2}} q_{bd};0_{n\times 1}]$ for Dirichlet boundary conditions. Here, $q_{bd}:=[q_0; 0_{(n-2)\times 1}; q_{n+1}]$ denotes the boundary term.  Time discretization of (\ref{expansion}) can be achieved by using the implicit symplectic integrator scheme (\ref{sischeme}). If $g(q)=0$, the successive over relaxation can be used to update the linear system for each time step; otherwise, the system is nonlinear and the Newton iteration can be used to time advance one step.

\subsection{Linear Wave Equation}
For our numerical experiments, we first study a linear system with $G(u)=g(u)=0$ and  with periodic boundary conditions. Let  $s(x)=10\times |x-{1\over 2} |$; and let $h(s)$ be a cubic spline function, which is  $1-\frac{3}{2}s^2+\frac{3}{4}s^3$ if $0\le s \le 1 $, $\frac{1}{4}(2-s)^3$ if $1<s\le2$, and 0 if $s>2$.
The initial condition is provided by
\begin{equation}\label{initial}
q(0)=[h(s(x_1));\ldots; h(s(x_n))], \quad p(0)=0_{n\times 1},
\end{equation}
which gives rise to a periodic system with wave propagating in both directions of $x$ in a periodic domain. The full model (reference benchmark solver) is computed using the following parameter set:
\medskip
\begin{center}
\begin{tabular}{r||l}
  \hline
  Space interval & $l=1$ \\
  Number of grid points & $n=500$ \\
  Space discretization step & $\Delta x=l/n=0.002$ \\
  Time interval & $T=50$\\
  Time discretization step & $\delta t=0.01$ \\
  Speed of the wave & $c=0.1$ \\
  \hline
\end{tabular}
\end{center}
\medskip

Regarding reduced systems, we compare all the proposed PSD algorithms (cotangent lift, complex SVD, and NLP with $r=100$) with the tangent lift in~\cite{Marsden:03h}, as well as with the POD. Both the Hamiltonian approach (PSD algorithms) and Lagrangian approach (tangent lift) are geometric algorithms, as they preserve the Hamiltonian or Lagrangian structures. The criterion for the comparison is the $L^2$ error between the generalized coordinate $q(t)$ of the benchmark solution and its approximations computed by reduced systems. For the Hamiltonian approach, PSD algorithms are used to construct symplectic matrices to fit the weighted data ensemble (\ref{gamma}) with $\gamma=\delta t$. For the Lagrangian approach,
 the basis matrix is constructed from a data ensemble of $q(t)$ in the configuration space. Thus, PSD reduced systems and the reduced system constructed by cotangent lift live on different subspaces of $\mathbb{R}^{2n}$. However, because $\dot q(t)= p(t)$ holds for the wave equation, symplectic integrators of $[q(t); p(t)]$ can also be used for the time  integration of $[q(t); \dot q(t)]$ for the fully-resolved Lagrangian system, and the reduced Lagrangian system constructed by the tangent lift also has the form (\ref{expansion}) if we assume $y(t)=[q(t); \dot q(t)]$.

\begin{figure}
\begin{center}
\begin{minipage}{0.48\linewidth} \begin{center}
\includegraphics[width=1\linewidth]{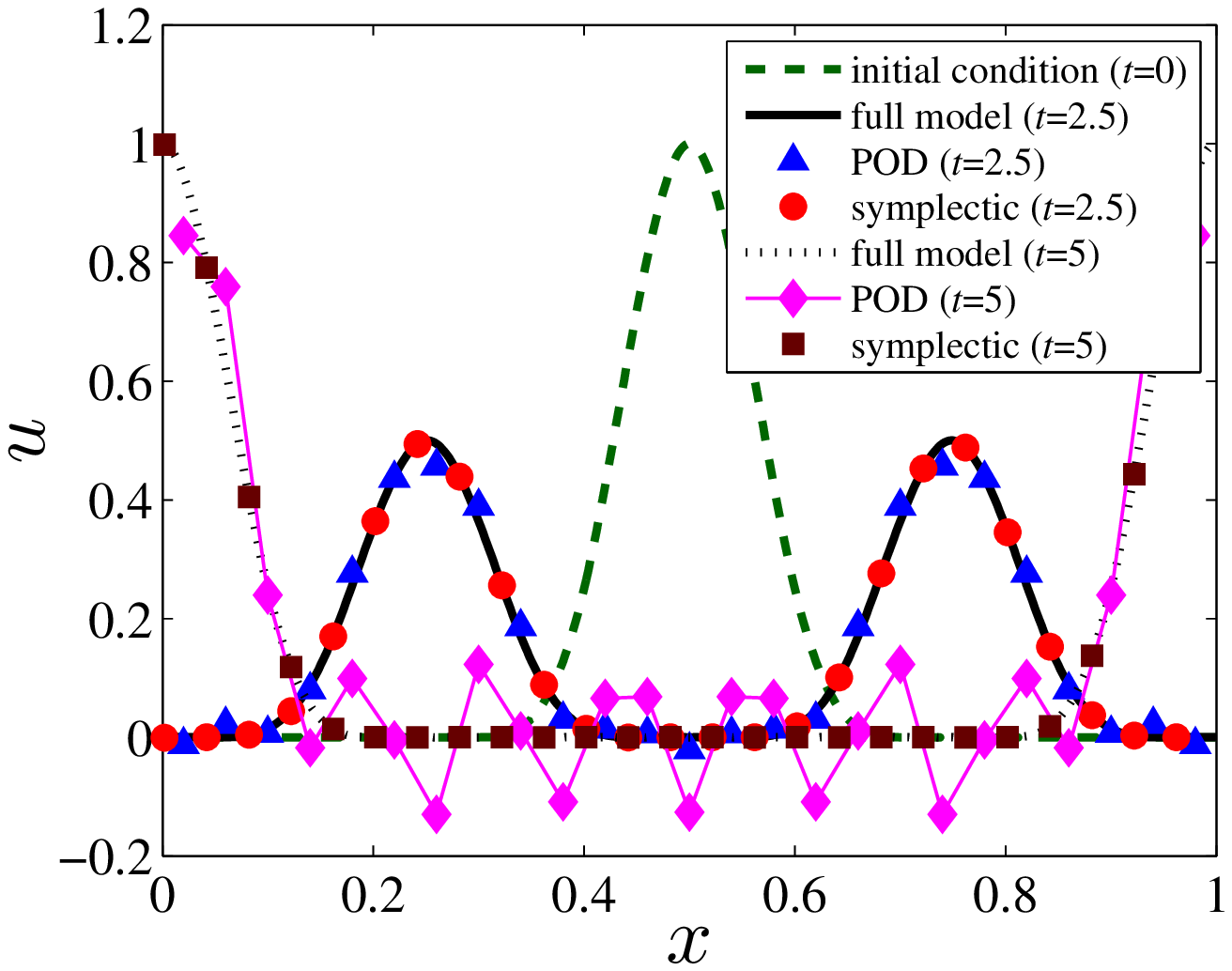}
\end{center} \end{minipage}
\begin{minipage}{0.48\linewidth} \begin{center}
\includegraphics[width=1\linewidth]{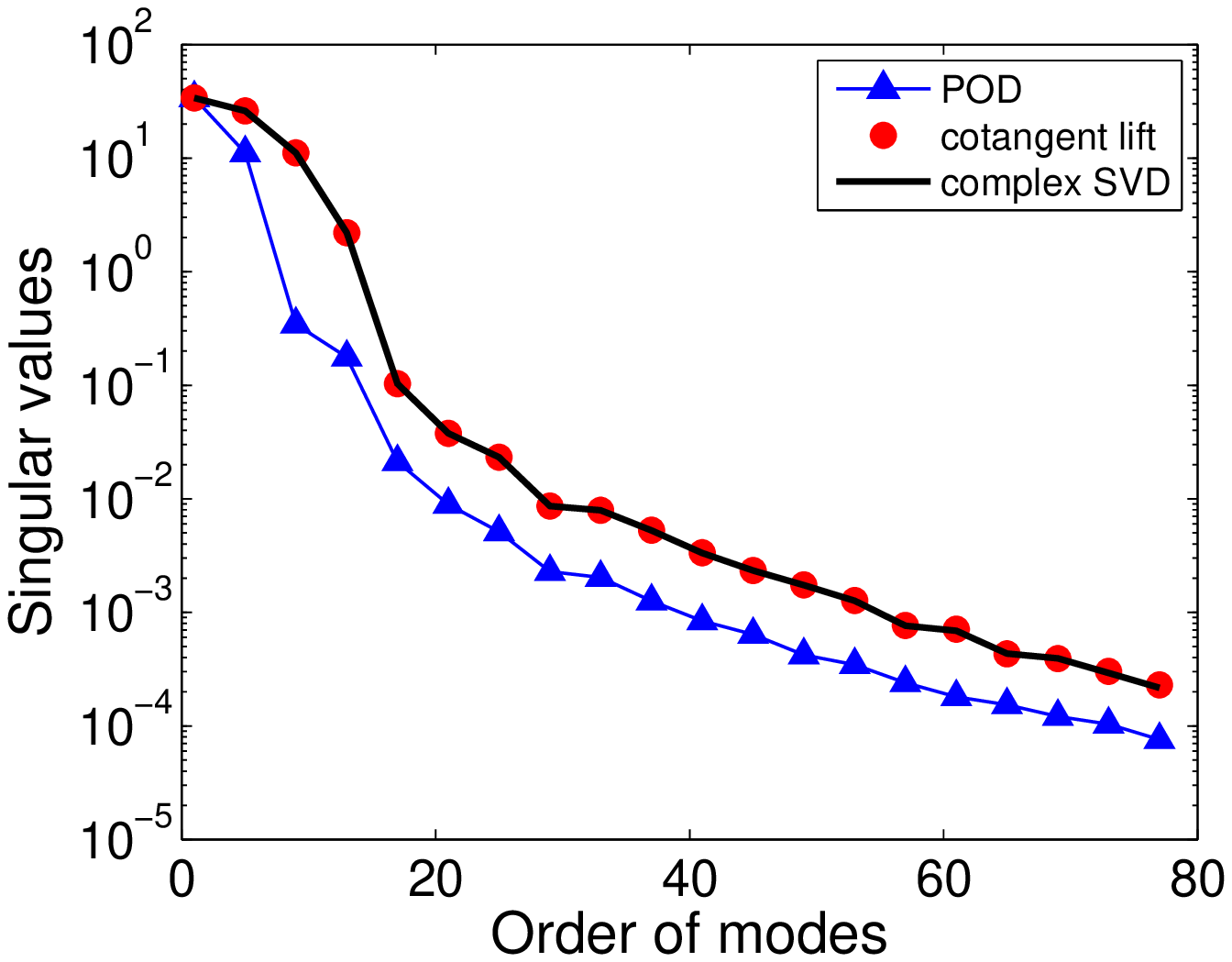}
\end{center} \end{minipage}\\
\begin{minipage}{0.48\linewidth}\begin{center} (a) \end{center}\end{minipage}
\begin{minipage}{0.48\linewidth}\begin{center} (b) \end{center}\end{minipage}
 \caption{(Color online.) (a) The solution $u(t,x)$ at $t=0$, $t=2.5$, and $t=5$ of the linear wave equation. We plot the results from the full model, the POD, and the cotangent lift.  (b) Plot the first 80 singular values corresponding to the first 80 POD (or PSD) modes that are used in different reduced systems.
} \label{fig:wavepod}
 \vspace{-3mm}
\end{center}
\end{figure}

Figure \ref{fig:wavepod}(a) shows the solution profile at $t = 0$, $t = 2.5$, and $t = 5$. The empirical data ensemble takes $101$ snapshots from the benchmark solution
trajectory with uniform interval ($\Delta t=0.5$). We first compare the POD with the cotangent lift.
 When $t=2.5$, both approaches can obtain good results by taking the first 20 modes; but when $t=5$, the POD reduced system significantly deviates from the full model.

In Figure \ref{fig:wavepod}(b), the blue line  represents the singular values  of the snapshot matrix $M_{x_\gamma}$ for the POD. Suppose $\{ \lambda_1, \ldots, \lambda_k\}$ denote the singular values of the snapshot matrix $M_1$ (or $M_2$) of the cotangent lift (or complex SVD). The red (or black) line represents the singular values $\{ \lambda_1, \lambda_1,  \ldots, \lambda_k, \lambda_k\}$ corresponding to the symplectic basis matrix $A_1={\rm{diag}}(\Phi, \Phi)$ (or $A_2=[\Phi, \Psi; -\Psi, \Phi]$). A fast decay of singular values indicates that a few modes can fit the data with good accuracy. This is a necessary (but not sufficient) condition for a low-dimensional reduced model to approximate the original system with good accuracy. Moreover, we notice that an arbitrary subspace of $\mathbb{R}^{2n}$ can be represented by an orthonormal basis matrix. However, unless this subspace is also symplectic, we cannot represent it by a symplectic basis matrix. Since PSD algorithms can only construct subspaces with the symplectic constraint, both the cotangent lift and complex SVD cannot fit the empirical data as well as the POD for the same subspace dimension.

\begin{figure}
\begin{center}
\begin{minipage}{0.48\linewidth} \begin{center}
\includegraphics[width=1\linewidth]{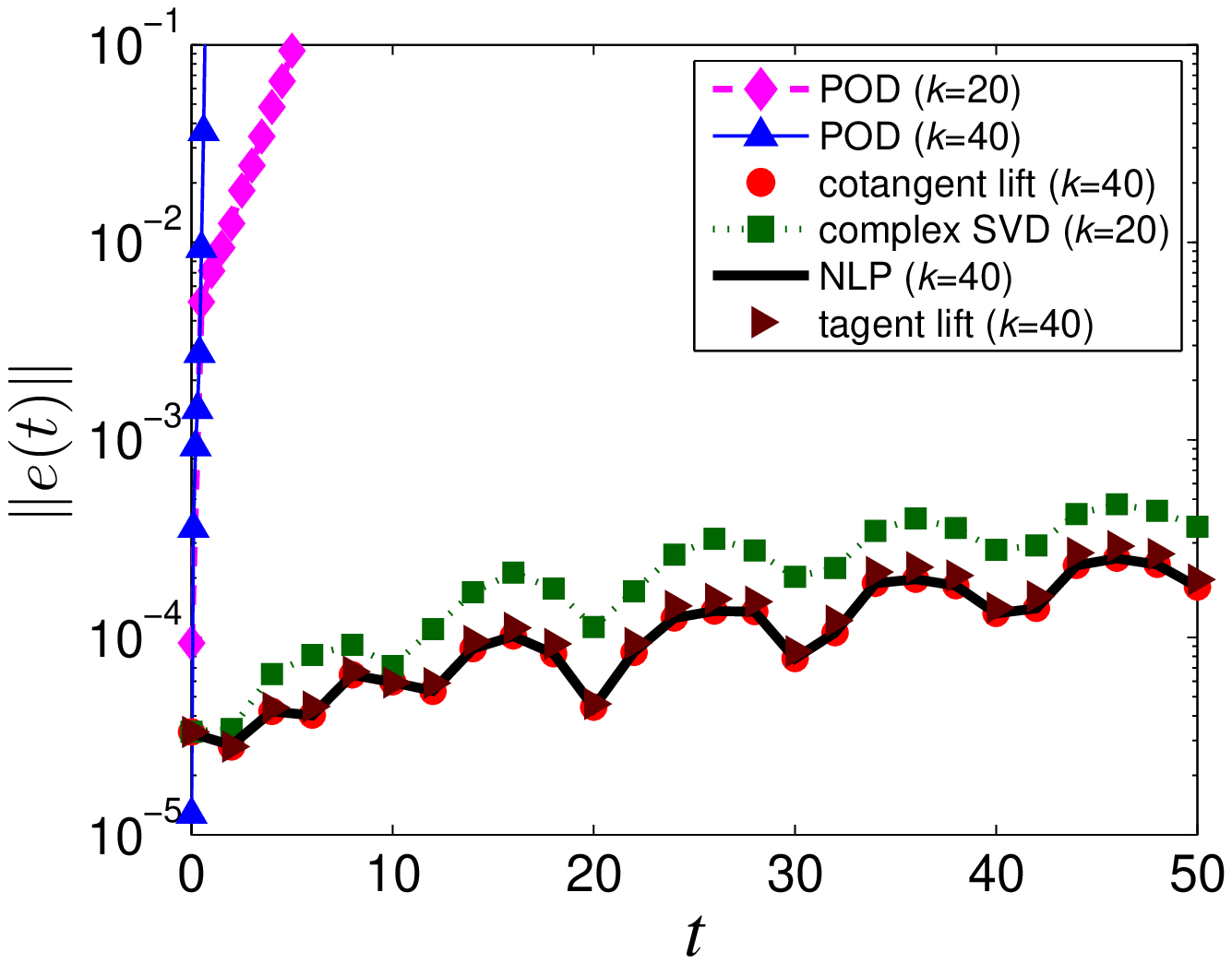}
\end{center} \end{minipage}
\begin{minipage}{0.48\linewidth} \begin{center}
\includegraphics[width=1\linewidth]{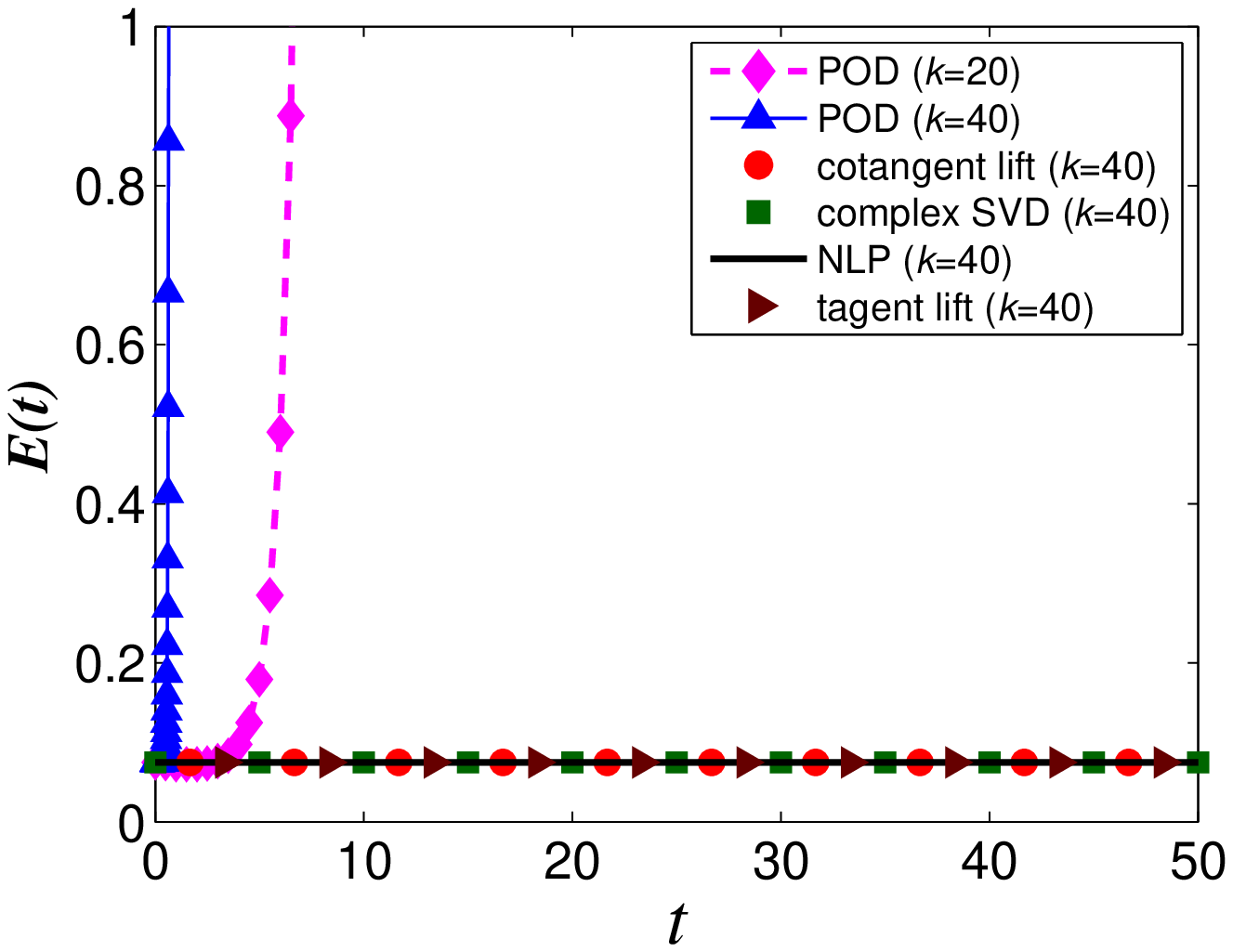}
\end{center} \end{minipage}\\
\begin{minipage}{0.48\linewidth}\begin{center} (a) \end{center}\end{minipage}
\begin{minipage}{0.48\linewidth}\begin{center} (b) \end{center}\end{minipage}
 \caption{(Color online.)   (a) The evolution of instant  $L^2$ error, $\|e(t)\|:=\|\hat u(t)-u(t)\|$, between the benchmark solution $u(t)$ and approximating solutions $\hat u(t)$  of the linear wave equation.  (b) The evolution of the energy $E(t)$ of different reduced systems.
} \label{fig:waveerr}
 \vspace{-3mm}
\end{center}
\end{figure}

\begin{figure}
\begin{center}
\begin{minipage}{0.5\linewidth} \begin{center}
\includegraphics[width=1\linewidth]{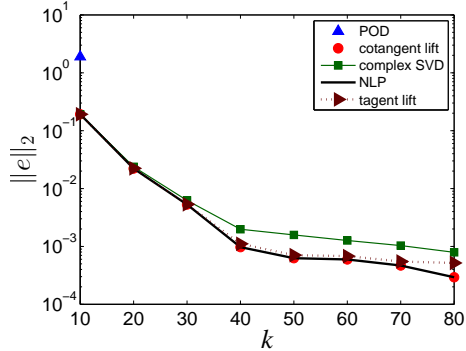}
\end{center} \end{minipage}\\
\caption{(Color online.) The $L^2$ norm of the total error $\|e\|_2:=\sqrt {\int_0^T {{{\left\| {e(t)} \right\|}^2}{\rm{d}}t} }$ of different reduced systems for the linear wave equation. For the POD reduced system, we only compute $\|e\|_2$ for $k=10$; when $k\ge20$, the reduced system blows up in the interested time domain [0, 50] and $\|e\|_2$ becomes infinite.} \label{fig:waveerr2}\vspace{-3mm}
\end{center}
\end{figure}

Using more modes, one may expect that both POD and PSD can produce more accurate solutions. However, as Figure \ref{fig:waveerr}(a) indicates, a POD reduced system blows up when it has 20 or 40 modes. In addition, the POD reduced system with 40 modes blows up faster than the system with 20 modes. This result verifies that the POD-Galerkin approach can yield unstable reduced systems, even though the original system is stable.  By contrast, errors in PSD reduced systems grow slowly in time. Figure \ref{fig:waveerr}(b) demonstrates that all the geometric algorithms preserve the system energy $E$, while the energy of POD reduced systems quickly grows to infinity. Here, $E$ equals the discretized Hamiltonian $H_d(y)$. Let $q_0=q_n$, (\ref{waveeng}) yields
\begin{equation}\label{hamlinear}
H_d(y) = \frac{\Delta x}{2} \sum\limits_{i = 1}^n p_i^2  +  \frac{c^2}{2\Delta x}\sum\limits_{i = 1}^n (q_i - q_{i - 1})^2.
\end{equation}

Figure \ref{fig:waveerr2} indicates that the $L^2$ norm of the total error of a POD reduced system is bounded in the interested time domain [0, 50] only for $k=10$ for the cases tested with $k\in \{10, 20, \ldots, 80\}$. While reduced systems constructed by  geometric algorithms show some numerical error, this error could be systematically reduced by using more modes. In terms of numerical accuracy, the cotangent lift and NLP are slightly better than the tangent lift of~\cite{Marsden:03h}, while the complex SVD is slightly worse than the tangent lift. The NLP yields the most accurate results, but for each $k$, we only observe a maximal of $0.028\%$ improvement compared with cotangent lift in terms of the relative percentage error.

\emph{Stability Preservation of Symplectic Model Reduction.} To explain our observations mentioned above, we study the stability of the linear wave equation. According to \cite{LeVequeRJ:07a}, the eigenvalues $\beta _i$ ($i= 1, \ldots, n$) of the discretized spatial derivative ${D}_{xx}$ with periodic boundary conditions are given by \[ \beta _i = -\frac{2}{\Delta x^2} \left[   1- \cos \left( \frac{2\pi i}{n} \right) \right],\] and the corresponding eigenvectors are given by \[{w_i} = \frac{1}{{\sqrt n }}\left[ e^{ - 2\pi \iota i/n},...,e^{ - 2\pi \iota i(n - 1)/n}, 1 \right].\] It follows that the eigenvalues of the Hamiltonian matrix $K$ in (\ref{wavematrix}) are given by $2n$ pure imaginary numbers $\pm \{\iota \gamma_i\}_{i=1}^n$, where $\gamma_i= c \sqrt{-\beta_i}$; and the corresponding eigenvectors are given by
$$
\xi _i:=\frac{1}{\sqrt{1+\gamma_i^2}}\begin{bmatrix}
 w_i\\
 \iota \gamma_i w_i
\end{bmatrix};\qquad
\zeta _i:=\frac{1}{\sqrt{1+\gamma_i^2}}
\begin{bmatrix}
 w_i\\
- \iota \gamma_i w_i
\end{bmatrix}.
$$
Since $\xi_n=\zeta_n={1 \over \sqrt{n}}[1_{n\times 1}; 0_{n\times 1}]$ by the above definition, we can redefine $\zeta_n$ to be $\zeta_n={1 \over \sqrt{n}}[0_{n\times 1}; 1_{n\times 1}]$. Thus, we can construct an invertible matrix $Q:=[\xi_1, \zeta_1, \ldots, \xi_n, \zeta_n ]$  such that $K$ is transformed to a Jordan form
\[Q^{-1}KQ={\rm{diag}}\left\{ \iota {\gamma_1}, -\iota {\gamma_1}, \dots ,\iota {\gamma_{n-1}}, - \iota {\gamma_{n - 1}}, \begin{bmatrix}
 0 & 1\\
0 & 0
\end{bmatrix} \right\}.\]
Although $K$ contains an unstable mode $\zeta_n$, the projection coefficient of initial condition (\ref{initial}) onto this mode vanishes, i.e., $\zeta_n^Ty_0=0$. Thus, the original system is stable and bounded for all $t$.

Next, we consider the reduced system constructed by the symplectic projection. By (\ref{hamlinear}), we have $H_d(y)\ge 0$, and the equality holds if and only if $y$ is parallel to $\xi_n$. If $\xi_n \notin {\rm{Range}}(A)$, the origin of $\mathbb{R}^{2n}$ is the only solution that satisfies $H_d(y)=0$ for all $y \in {\rm{Range}}(A)$, as a result, it is a strict minimum of $H_d$.  In our numerical simulations, we do observe that $A^+\xi_n\ne 0$, where $A$ is constructed by the cotangent lift or complex SVD. Thus, by the Dirichlet's stability theorem, the origin is a stable equilibrium for the reduced Hamiltonian system, which implies that the symplectic projection preserves the stability of the linear wave equation.

\emph{Instability of POD-Galerkin.}
Since the POD does not preserve the system energy, there are no mechanisms similar to the Hamiltonian and Lagrangian approaches that limit the solution trajectory in a bounded region.  As a result, the reduced system may blow up with time evolution.  To corroborate this claim, let $\lambda_*$ denote the eigenvalue of $\Phi^TK \Phi$ with the maximal real part and let $\xi_*$ denote the corresponding eigenvector with unit length.  Then, $a_*=\xi_*^T y_0$ gives the projection coefficient of $y_0$ onto $\xi_*$. The following table indicates that for different subspace dimensions $k$, a POD reduced system has ${\rm{Re}}(\lambda_*)>0$ and $a_*\ne 0$.  Since the solution of a linear system has an exponential term $a_* \exp(\lambda_* t)\xi_*$, the POD reduced system is always unstable for long-time integration.

\begin{table} [htbp]
\begin{center}
 \label{tab:podeigv}
\begin{tabular}{|r|c|c|c|c|c|c|c|c|}
\hline
 $k$ &  10 & 20 &30 &40&50&60&70 &80  \\
\hline
   $\lambda_*$ & 0.0338 & 0.659 & 13.74 & 14.39 & 14.50 & 5.33 & 10.42 & 13.05 + 5.09 $\iota$\\
  \hline
$a_*$ &   0.929& 0.0184 & 0.0263& 0.0498& -0.0200& 0.0718 & 7.55e-3 & -0.0068 - 0.0142$\iota$\\
  \hline
\end{tabular}
\end{center}
\end{table}

Assume that $\Phi_{k_1}$ and $\Phi_{k_2}$ respectively contain the first $k_1$ and $k_2$ dominant modes. If $k_1<k_2$, then $\Phi_{k_1}^T K \Phi_{k_1}$ is a submatrix of $\Phi_{k_2}^T K \Phi_{k_2}$, and $\|\Phi_{k_1}^T K \Phi_{k_1}\| \le \|\Phi_{k_2}^T K \Phi_{k_2}\|$ holds. As ${\rm{Re}}(\lambda_*)\le |\lambda_*| \le \|\Phi^T K \Phi\|$, the matrix norm of $\Phi^T K \Phi$ provides an upper bound for ${\rm{Re}}(\lambda_*)$. Thus, the upper bound of ${\rm{Re}}(\lambda_*)$  is a monotonically increasing function of the subspace dimension $k$. The above table also shows that ${\rm{Re}}(\lambda_*)$ with 40 modes is much larger than ${\rm{Re}}(\lambda_*)$ with 20 modes, which explains why the POD reduced system with 40 modes blows up faster than the system with 20 modes in Figure \ref{fig:waveerr}.  Although for $k=10$, POD can produce a reduced system with reasonable accuracy for a short time domain [0, 2.5], we can still observe that for a large enough integration time, say $t>10$, this system blows up.

\subsection{Sine-Gordon Equation}
Next, we consider a special nonlinear wave equation with $G(u)=1-\cos(u)$, $g(u)=\sin (u)$ and $c=1$, which corresponds to the sine-Gordon equation.  This equation, which was first studied  in the 1970s, appears in a number of physical applications, including relativistic field theory, Josephson junctions, and mechanical transmission lines~\cite{WhithamGB:99a}. One can show that the sine-Gordon equation admits a localized solitary wave solution,
\begin{equation}\label{gsolution}\Scale[0.95]{
u(t,x) = 4\arctan \left[ {\exp \left( { \pm \frac{{x - {x_0} - vt}}{{\sqrt {1 - {v^2}} }}} \right)} \right]},
\end{equation}
which travels with the speed $|v|<1$. The $\pm$ signs correspond to localized solutions which are called \emph{kink} and \emph{antikink}, respectively~\cite{WhithamGB:99a}.

In our simulations, the full model is solved for the kink case with Dirichlet boundary conditions ($u(t,0)=0, u(t, 1)=2 \pi$) using the following parameter set:
\medskip
\begin{center}
\begin{tabular}{r||l}
  \hline
  Space interval & $l=50$ \\
  Number of grid points & $n=2000$ \\
  Space discretization step & $\Delta x=l/n=0.025$ \\
  Time interval & $T=150$\\
  Time discretization step & $\delta t=0.0125$ \\
  Speed of the wave & $v=0.2$ \\
  \hline
\end{tabular}
\end{center}
\medskip

 The $L^2$ error for the state variable $y(t)$  is studied for the full model and reduced models constructed by the POD, cotangent lift, complex SVD,  DEIM, and  SDEIM.
All the basis matrices are constructed to fit the data ensemble (\ref{gamma}) with $\gamma=1$.  For the SDEIM, the cotangent lift is used to construct a symplectic basis matrix $A_1=\Psi={\rm{diag}}(\Phi, \Phi)$ to fit both the state variable $y(t)$ and nonlinear vector $f_N(y(t))$, where $\Phi$ is the POD basis matrix for the extended snapshot matrix that contains $q(t)$, $p(t)$, and $f_N(q(t))$ in its column vectors.

\begin{figure}
\begin{center}
\begin{minipage}{0.48\linewidth} \begin{center}
\includegraphics[width=1\linewidth]{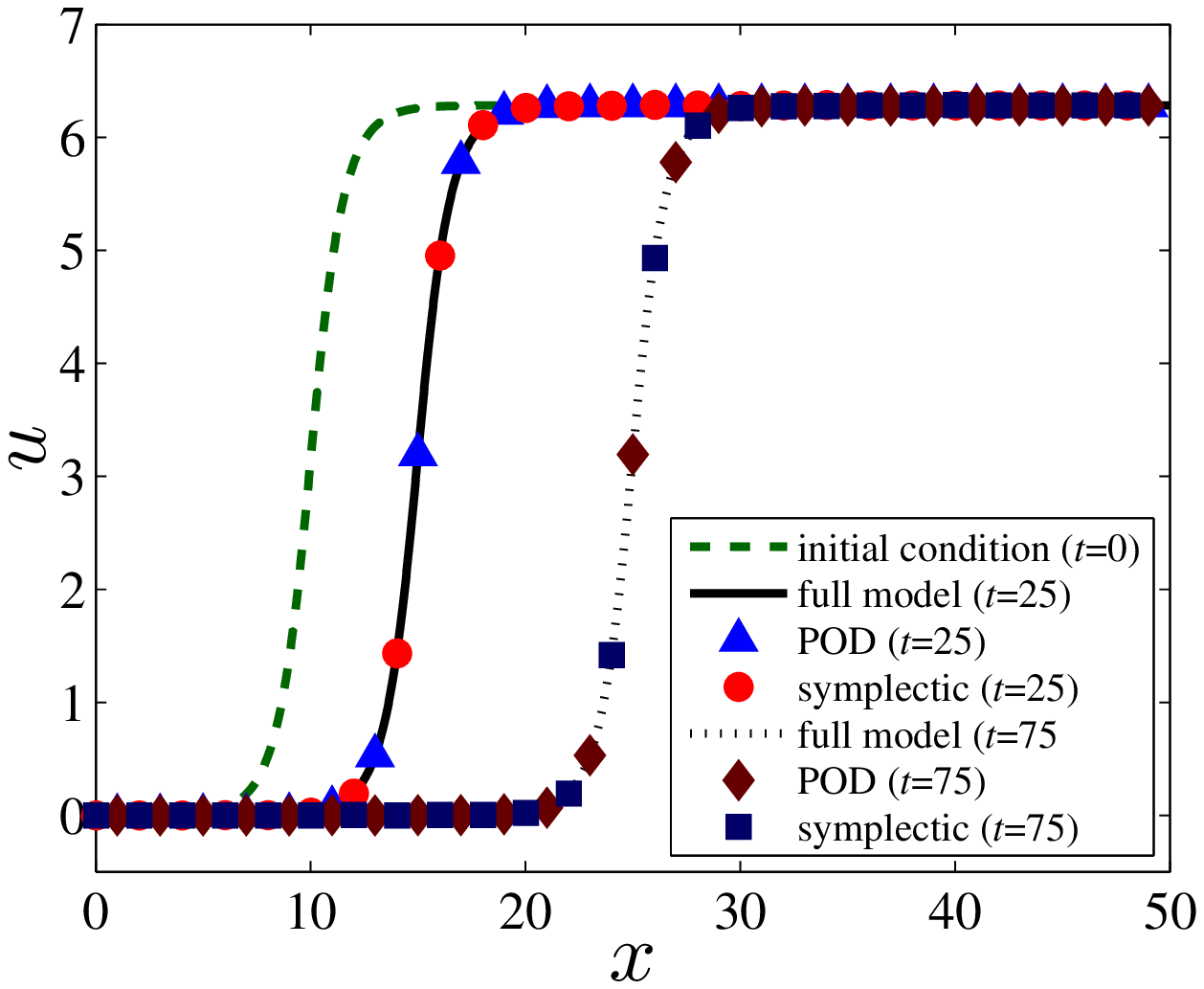}
\end{center} \end{minipage}
\begin{minipage}{0.48\linewidth} \begin{center}
\includegraphics[width=1\linewidth]{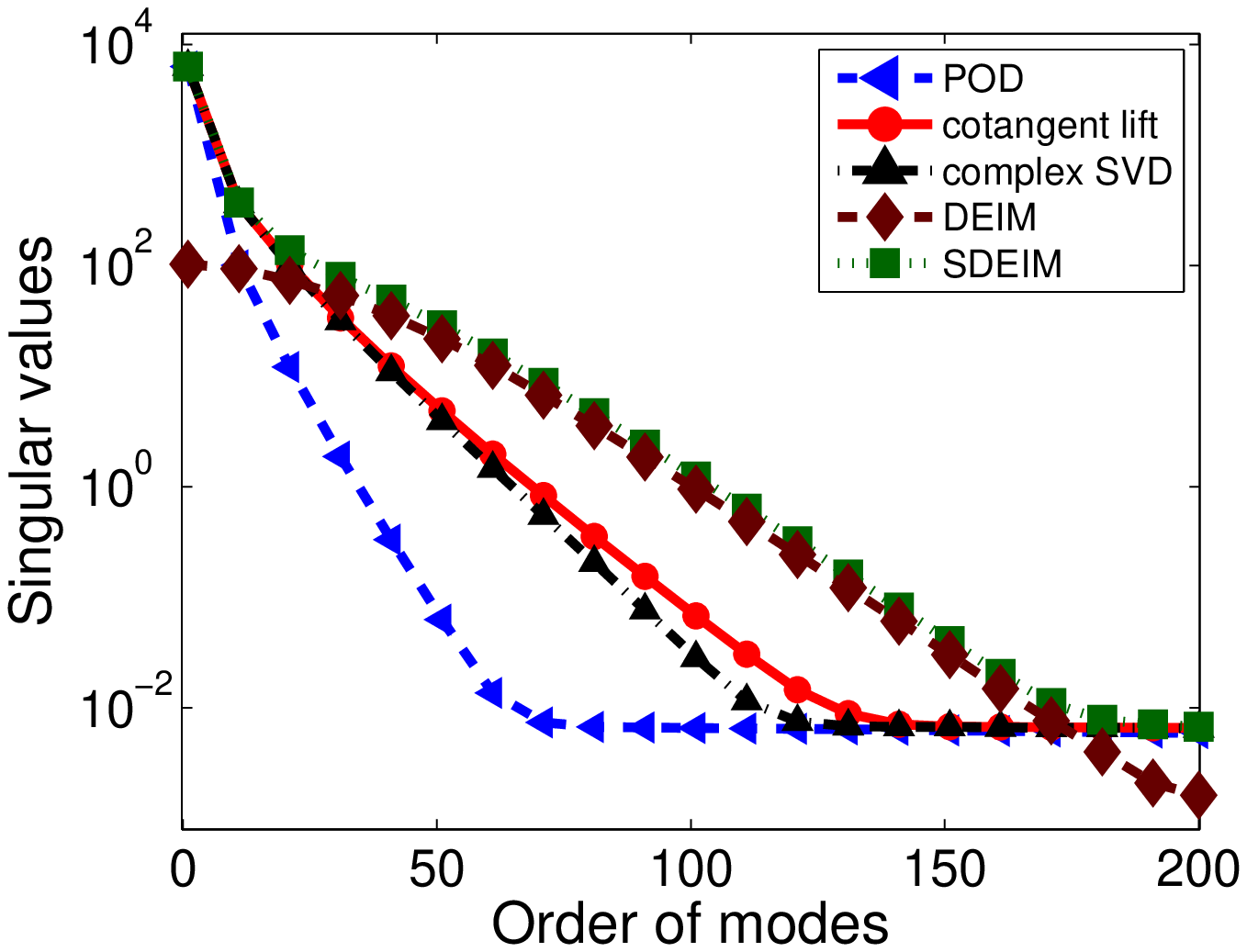}
\end{center} \end{minipage}\\
\begin{minipage}{0.48\linewidth}\begin{center} (a) \end{center}\end{minipage}
\begin{minipage}{0.48\linewidth}\begin{center} (b) \end{center}\end{minipage}
 \caption{(Color online.) (a) The solution $u(t,x)$ at $t=0$, $t=25$, and
$t=75$ of the sine-Gordon equation.  We plot the results from the full model, the POD, and the cotangent lift. (b) Plot the first 200 singular values corresponding to the first 200 POD (or DEIM, PSD, SDEIM) modes that are used in different reduced systems.
} \label{fig:gpod}
 \vspace{-3mm}
\end{center}
\end{figure}

Figure \ref{fig:gpod}(a) shows the \emph{kink} solution profile at $t=0$, $t = 25$, and $t=75$. The data ensemble takes $1201$ snapshots from the solution trajectory, solved by the full model with uniform interval ($\Delta t=0.125$). We first compare the POD with the cotangent lift.  For short-time integration, both approaches could obtain very accurate results by taking the first 60 modes. In Figure \ref{fig:gpod}(b), we study the  singular values corresponding to the POD basis matrix, the PSD basis matrices constructed by the cotangent lift and complex SVD,  the nonlinear term basis matrix  for the DEIM, and the symplectic matrix for the SDEIM.  This figure demonstrates that the POD is better to fit empirical state variables than the cotangent lift and complex SVD, while the DEIM is better to fit empirical nonlinear vectors than the SDEIM.

\begin{figure}
\begin{center}
\begin{minipage}{0.48\linewidth} \begin{center}
\includegraphics[width=1\linewidth]{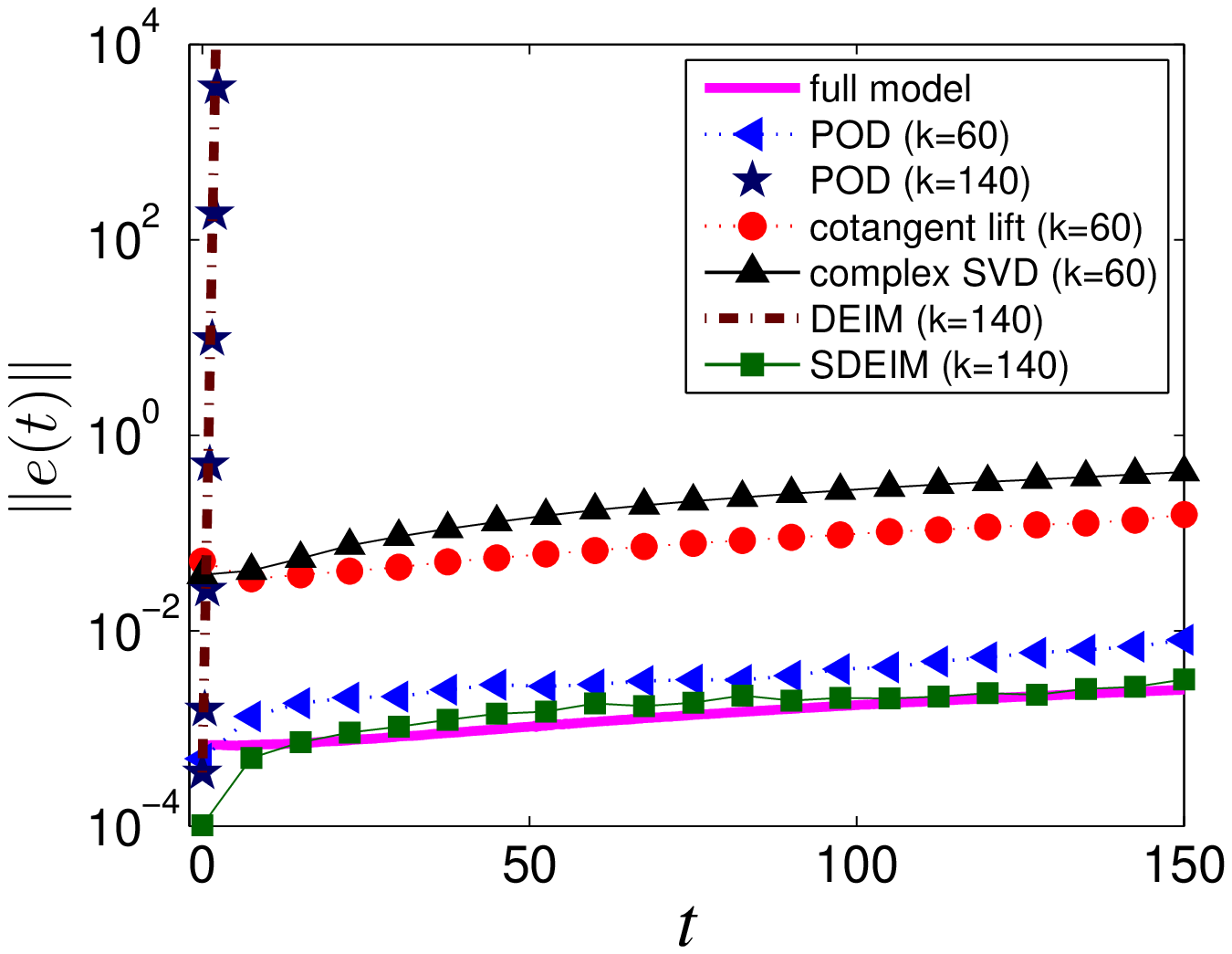}
\end{center} \end{minipage}
\begin{minipage}{0.48\linewidth} \begin{center}
\includegraphics[width=1\linewidth]{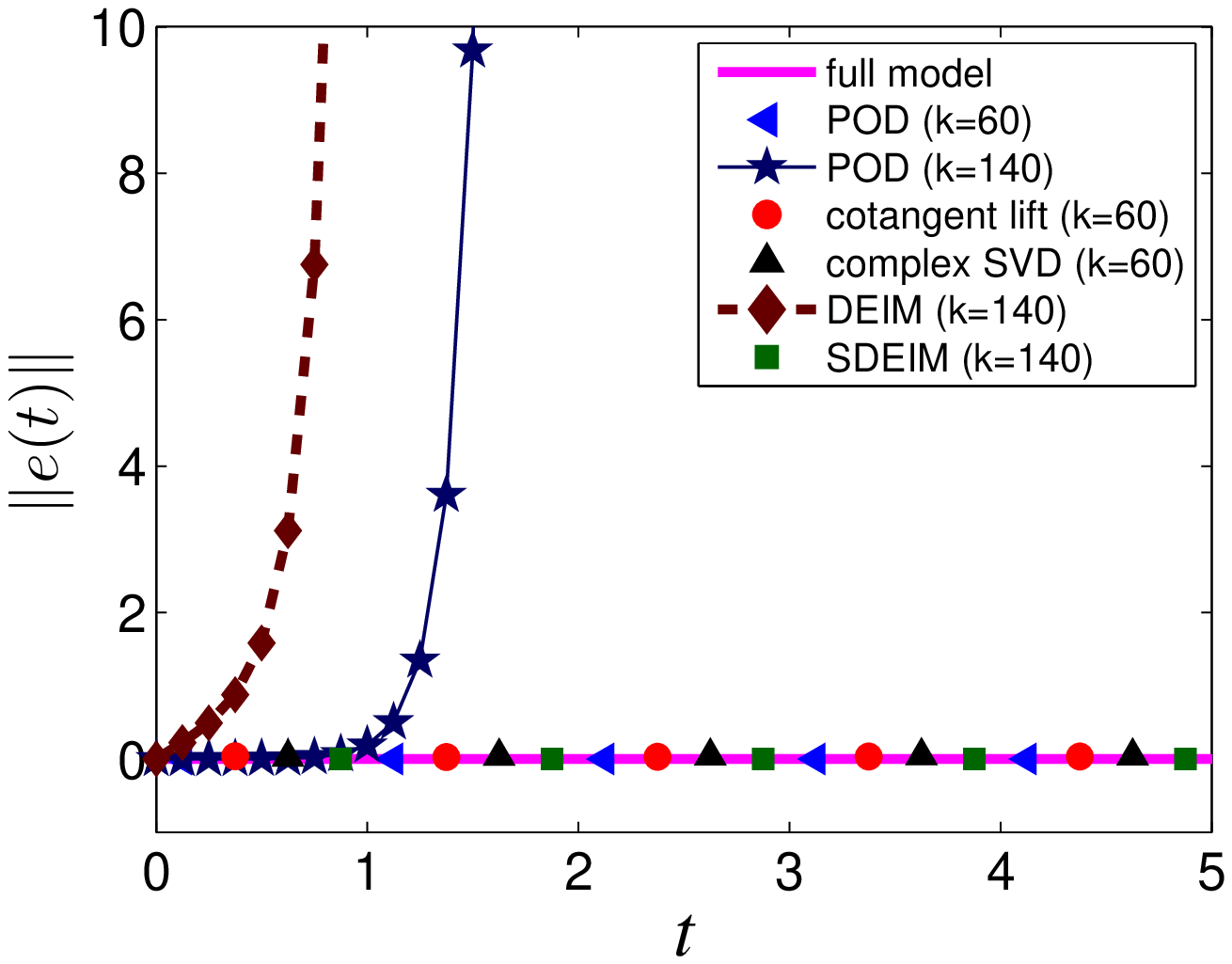}
\end{center} \end{minipage}\\
\begin{minipage}{0.48\linewidth}\begin{center} (a) \end{center}\end{minipage}
\begin{minipage}{0.48\linewidth}\begin{center} (b) \end{center}\end{minipage}
 \caption{(Color online.)   (a) The evolution of instant  $L^2$ error, $\|e(t)\|:=\|\hat
y(t)-y(t)\|$, between  the analytic solution $y(t)$ in the phase space and approximating solutions $\hat y(t)$  of the sine-Gordon  equation for $t\in[0, 150]$. (b) The instant $L^2$ error $\|e(t)\|$ for $t\in[0, 5]$.
} \label{fig:gerr}
 \vspace{-3mm}
\end{center}
\end{figure}

\begin{figure}
\begin{center}
\begin{minipage}{0.48\linewidth} \begin{center}
\includegraphics[width=1\linewidth]{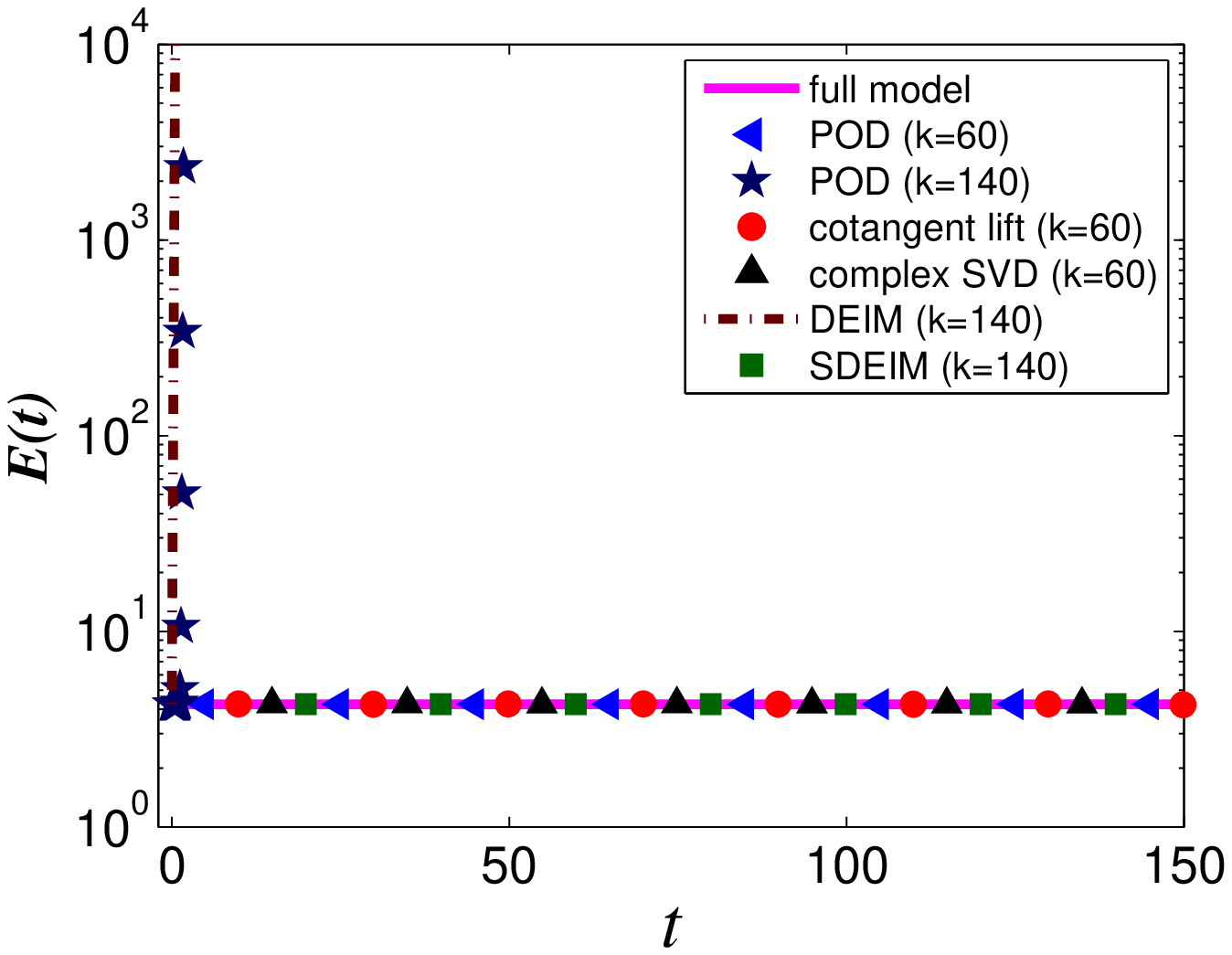}
\end{center} \end{minipage}
\begin{minipage}{0.48\linewidth} \begin{center}
\includegraphics[width=1\linewidth]{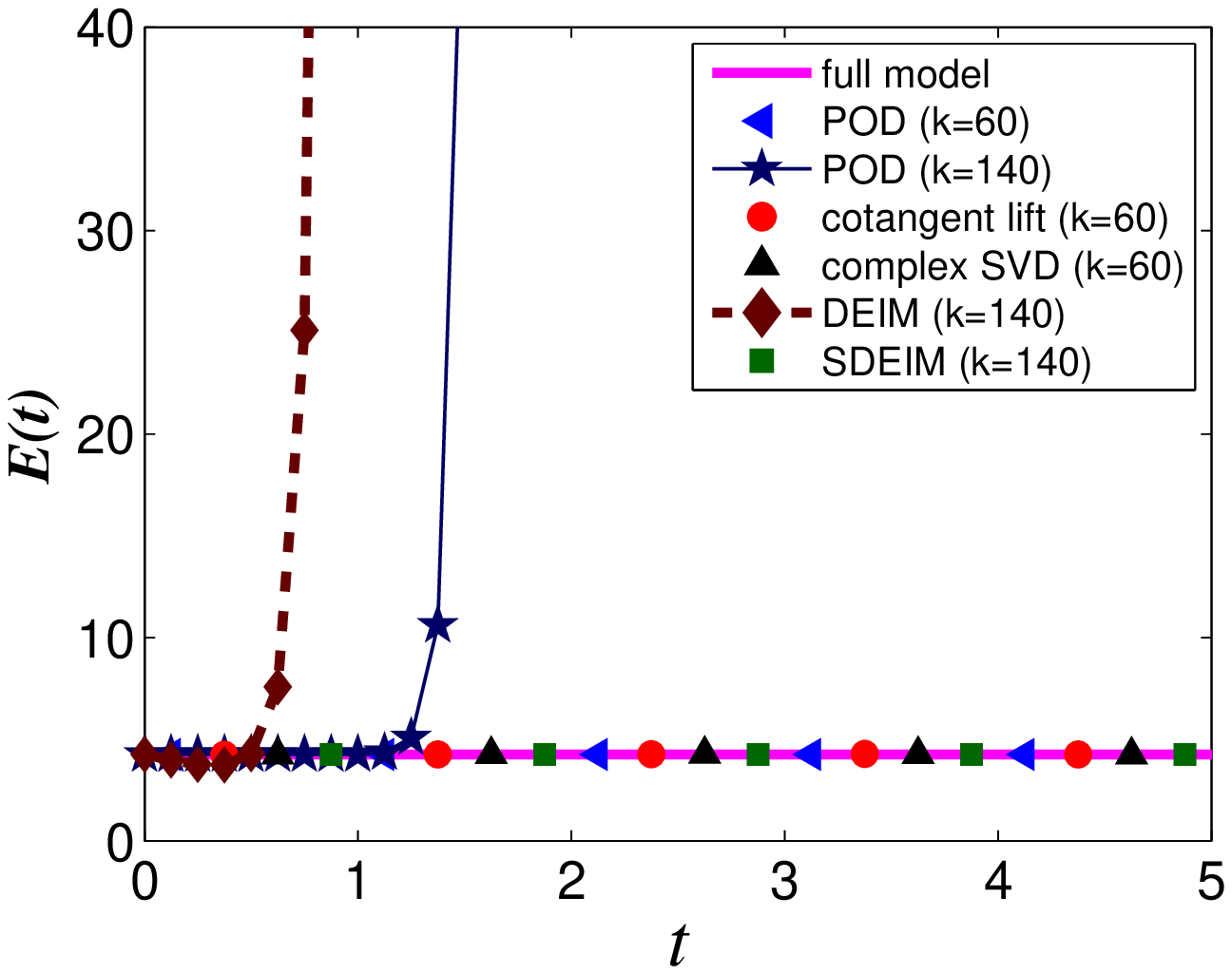}
\end{center} \end{minipage}\\
\begin{minipage}{0.48\linewidth}\begin{center} (a) \end{center}\end{minipage}
\begin{minipage}{0.48\linewidth}\begin{center} (b) \end{center}\end{minipage}
 \caption{(Color online.)
 (a) The evolution of the system energy $E(t)$ of the sine-Gordon equation for $t\in[0, 150]$. (b) The evolution of the system energy $E(t)$  for $t\in[0, 5]$.
} \label{fig:geng}
 \vspace{-3mm}
\end{center}
\end{figure}

Figure \ref{fig:gerr}  illustrates that all  symplectic schemes (including the cotangent lift,  complex SVD, and SDEIM) yield low computational errors with appropriate subspace dimension,
 while non-symplectic schemes (including the  POD and  DEIM) can yield unbounded numerical error with 140 modes. In Figure \ref{fig:geng},
 all symplectic schemes can effectively preserve the system energy $E$.
  By contrast, both POD and DEIM reduced systems can achieve infinite energy with 140 modes. Here, $E$ equals the discretized Hamiltonian $H_d(y)$. With $G(u)=1-\cos(u)$ and Dirichlet boundary conditions,  (\ref{waveeng}) gives
\begin{equation}\label{hamsin}
\Scale[0.96]{
H_d(y) =  \frac{\Delta x}{2} \sum\limits_{i = 1}^np_i^2  + \Delta x\sum\limits_{i = 1}^n [1 - \cos (q_i) ] + \frac{q_1^2}{4\Delta x} + \frac{1}{2\Delta x}\sum\limits_{i = 2}^n (q_i - q_{i - 1})^2  + \frac{(q_n - 2\pi )^2}{4\Delta x}}.
\end{equation}

\begin{figure}
\begin{center}
\begin{minipage}{0.48\linewidth} \begin{center}
\includegraphics[width=1\linewidth]{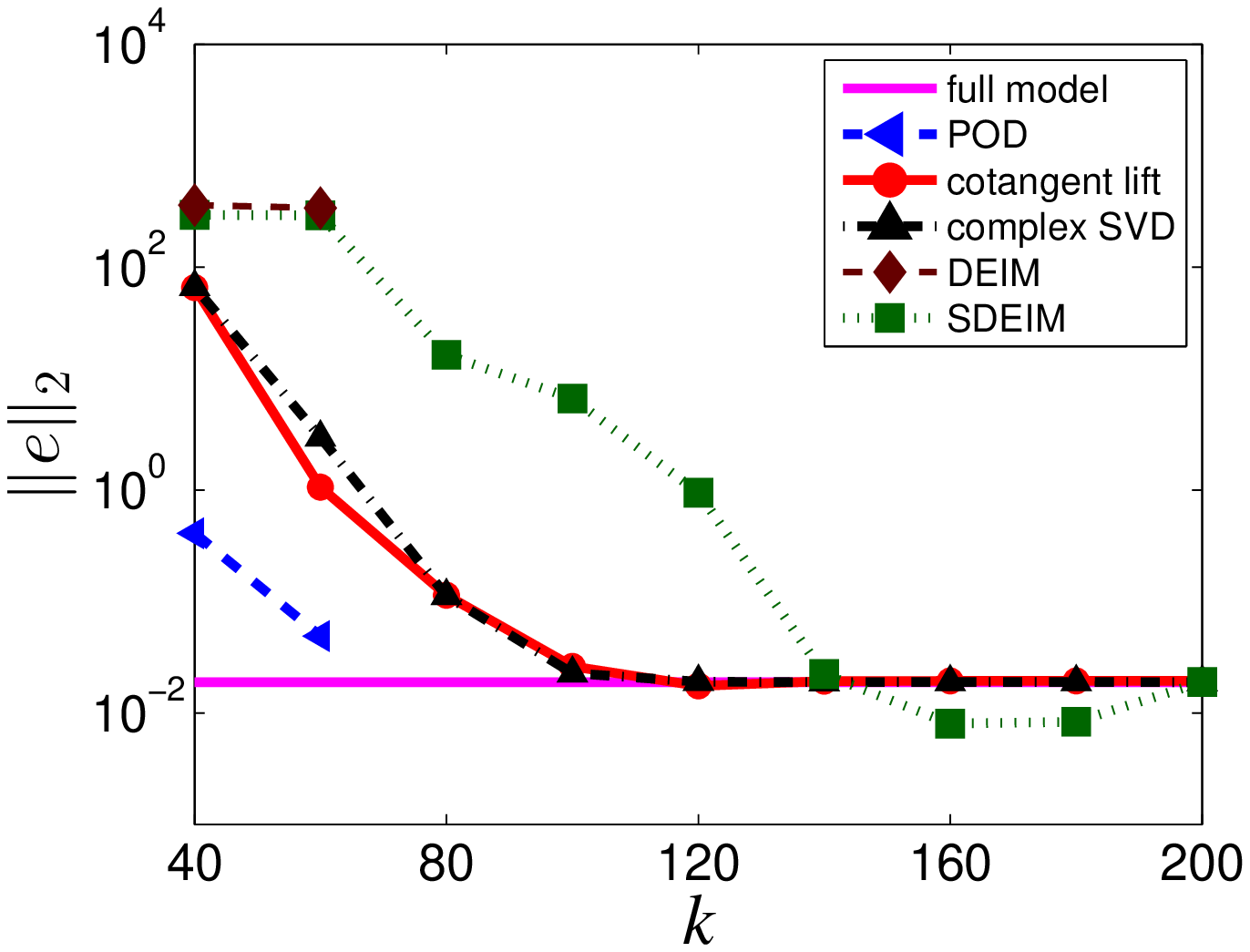}
\end{center} \end{minipage}
\begin{minipage}{0.48\linewidth} \begin{center}
\includegraphics[width=1\linewidth]{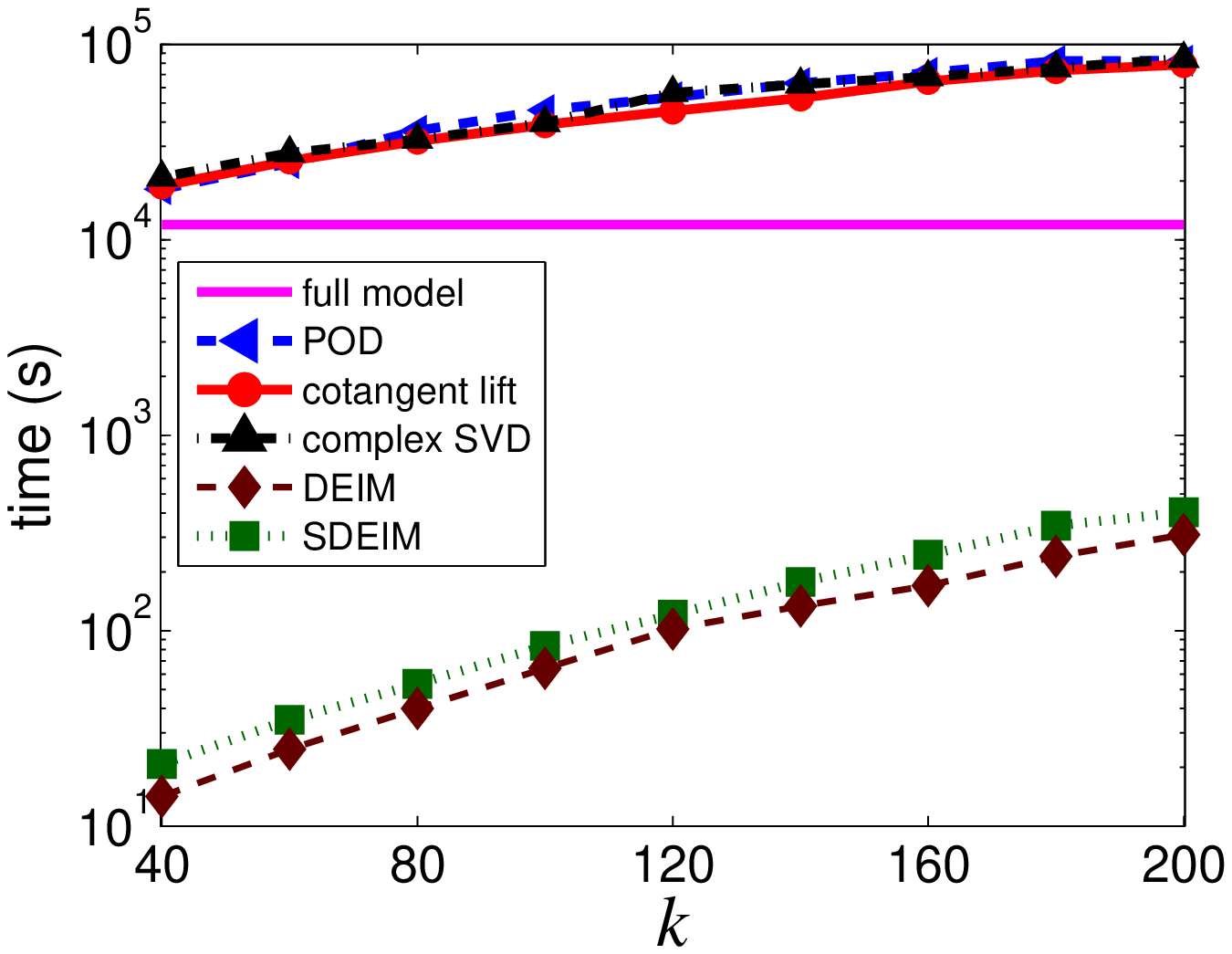}
\end{center} \end{minipage}\\
\begin{minipage}{0.48\linewidth}\begin{center} (a) \end{center}\end{minipage}
\begin{minipage}{0.48\linewidth}\begin{center} (b) \end{center}\end{minipage}
 \caption{(Color online.)
  (a) The $L^2$ norm of the total error $\|e\|_2:=\sqrt {\int_0^T {{{\left\| {e(t)} \right\|}^2}{\rm{d}}t} }$ of the full model and different reduced models for the sine-Gordon equation.  For the POD and DEIM reduced systems, we only compute $\|e\|_2$ for $k=40$ and $k=60$; when $k\ge 80$, the reduced systems blow up in the interested time domain [0,150] and $\|e(t)\|_2$ becomes infinite. (b) The running time of different model reduction techniques with different $k$. All the data come from the average value of ten independent runs.
} \label{fig:gtime}
 \vspace{-3mm}
\end{center}
\end{figure}

Figure \ref{fig:gtime}(a) indicates that by using more modes, all symplectic reduced models can obtain better accuracy and finally converge to the full model. By contrast, the POD and DEIM can yield unbounded reduced systems in the interested time domain [0, 150] for $k\ge 80$  for the cases tested with $k \in \{40,60,...,200\}$. Furthermore, by the analysis in Section \ref{sec:deim}, we know that a direct use of the POD-Galerkin or PSD with the symplectic projection is not able to obtain any speedups for the sine--Gordon equation, since it contains a nonlinear vector term $\sin(u)$. Numerical results in Figure \ref{fig:gtime}(b) also verify this point. Especially, the running time for the POD, the cotangent lift, and the complex SVD is even larger than the running time for the full model. On the other hand, both the DEIM and SDEIM approximations could significantly improve the computational
 efficiency and reduce the running time of POD or PSD by three orders of magnitude.

The boundedness of symplectic reduced systems can be derived by their energy conservation property. If $y(t)$ denote the solution trajectory, we have $H_d(y(t))=E$ for a constant $E$. Since each term on the RHS of (\ref{hamsin}) is nonnegative, we must have $|p_i|\le \sqrt{2E/\Delta x}$, $|q_1|\le 2\sqrt{E\Delta x }$, and $|q_i| \le |q_{i-1}|+\sqrt{2E\Delta x }$ for $i\ge 2$. In other words, there exists a positive number $M$, for any state $y\in \mathbb{R}^{2n}$, as long as $\|y\|=M$, we have $H_d(y)>E$. Therefore, by Theorem \ref{thm:bound} both the original system and  reduced systems constructed by the symplectic projection are bounded for all $t$.

\section{Conclusion}\label{sec:conclusion}
In this paper, we proposed a symplectic model reduction technique for the reduced-order modeling of large-scale Hamiltonian systems. We first defined the symplectic projection, which can yield reduced systems that remain Hamiltonian. Several proper symplectic decomposition (PSD) algorithms, such as the cotangent lift, complex SVD, and   nonlinear programming, were developed to generate a symplectic matrix that spans a low-dimensional symplectic subspace.

 Because the symplectic model reduction preserves the symplectic structure, it also preserves the system energy and stability. Thus, the proposed technique is very suitable for long-time integration, especially when the original systems are conservative and do not have any natural dissipative mechanism to stabilize them. Since the symplectic projection can only speed up linear and quadratic problems, the PSD was also combined with DEIM, effectively reducing the complexity of the nonlinear vector term.  Because the complexity of the symplectic discrete empirical interpolation method (SDEIM) does not depend on the dimension of the original system, a significant speedup can be obtained for a general nonlinear problem. We demonstrated the capability of the symplectic model reduction to solve a large-scale system with high accuracy, good efficiency, and stability preservation via linear and nonlinear wave equations.

\bibliographystyle{siam}
\bibliography{RefA4}

\end{document}